\documentclass[12pt,a4paper]{article}

\usepackage[ansinew]{inputenc}

\usepackage[french]{babel}

\usepackage{amssymb,amsmath, amsfonts}

\usepackage [dvips] {graphicx}

\usepackage{amsthm}

\usepackage{stmaryrd}

\newtheorem{df}{Définition}
\newtheorem{prop} {Proposition}
\newtheorem{lm} [prop]{Lemme}
\newtheorem{thm} [prop] {Théorème}
\newtheorem{cor} [prop] {Corollaire}

\theoremstyle{remark}
\newtheorem{rmq}{Remarque} 
\newtheorem{exm}{Exemple}

\author{Stéphane Dugowson}
\date{}

\title{Espaces connectifs : représentations, feuilletages, ordres et difféologies}

\begin{document}
 
\maketitle

Cet article est publié en juin 2013 dans les \textit{Cahiers de topologie et géométrie différentielle catégoriques}, volume LIV, avec cette dédicace :

\begin{flushright}
\textit{Dédié à René Guitart, en toute amitié}
\end{flushright}

\tableofcontents
\newpage



\vspace*{35mm}
\begin{small}
\begin{tabular}{@{} p{0.5cm} @{}  p{11.3cm} @{}}
& \textbf{Abstract}. This article is a continuation of my former article ``On Connectivity Spaces"\,\cite{Dugowson:201012}.  After some brief historical references relating to the subject, separation spaces and then adjoint notions of connective representation and connective foliation are developed. The connectivity order previously defined only in the finite case is now generalised to all connectivity spaces, and so to connective foliations. Finally, we start the study of some functorial relations between connectivity and diffeological spaces,  and we give a characterization of diffeologisable connectivity spaces. \\
& \textbf{Résumé.}  Poursuivant l'étude présentée dans notre article \og On connectivity Spaces\fg \,\cite{Dugowson:201012}, nous développons ici, après quelques rapides repères historiques, la notion d'espace de séparation et les notions adjointes de représentation et de feuilletage connectifs. Nous généralisons en outre la notion d'ordre connectif au cas infini, ainsi qu'aux feuilletages connectifs. Finalement, nous étudions certaines relations fonctorielles entre structures connectives et structures difféologiques, caractérisant en particulier les espaces connectifs difféologisables. \\
& \textbf{Keywords.} Connectivity. Links. Borromean. Foliation. Connective representation. Connectivity order. Diffeology\\
& \textbf{Mathematics Subject Classification (2010).}   54A05. 54B30. 57M25. 57R30. 58A99.
\end{tabular}
\end{small} 


\mbox{}


Les propriétés connectives du  \emph{n\oe ud borroméen} --- trois courbes globalement entrelacées mais libres deux à deux --- président, dès sa naissance en 1892 avec l'article fondateur du mathématicien allemand Hermann Brunn \cite{Brunn:1892a}, à l'histoire des espaces connectifs. Peu connu, le résultat de Hermann Brunn concernant la possibilité de représenter par entrelacs toute structure connective finie devra attendre les travaux du Suisse Hans Debrunner \cite{Debrunner:19600416,Debrunner:1964} dans les années 1960, puis ceux, vingt ans plus tard, du Japonais Taizo Kanenobu \cite{Kanenobu:198504,Kanenobu:1986} pour être rigoureusement établi, constituant ce que j'ai nommé le théorème de Brunn-Debrunner-Kanenobu\footnote{Voir \cite{Dugowson:201012}.}. Il faut dire que la notion d'espace connectif elle-même n'aura clairement été dégagée en tant que telle qu'au début des années 1980, par le mathémati\-cien allemand Reinhard Börger \cite{Borger:1981,Borger:1983} qui définit la catégorie des espaces connectifs\footnote{Plus précisément, Börger considère ce que nous appelons les espaces connectifs intègres, dans lesquels chaque point est nécessairement connexe.} et en donne les premières propriétés, en particulier le fait qu'il s'agit d'une catégorie complète et co-complète, mais sans faire le lien avec les travaux antérieurs de Brunn et Debrunner. Il semble, c'est du moins notre hypothèse, que l'intérêt de Börger pour cette catégorie ait été freiné par le constat qu'elle n'est pas cartésienne fermée\footnote{C'est par contre une catégorie monoïdale fermée, voir \cite{Dugowson:201012}.}. Par ailleurs, il se pourrait qu'en nommant \og entrelacs brunniens\fg\, les entrelacs effectivement utilisés par Brunn dans sa construction, mais en oubliant de signaler cette construction elle-même, Rolfsen ait, dans son ouvrage publié en 1976 sur la théorie des {n\oe uds}\cite{Rolfsen:1976}, rendu à Brunn un hommage un peu paradoxal, en ce sens que le souvenir des briques aura pu contribuer à l'oubli de la bâtisse. En 1988, dans le cadre de travaux sur l'analyse des images et la morphologie mathématique, les Français  Georges Matheron et Jean Serra \cite{Matheron.Serra:1988, Matheron.Serra:1988k}, ignorant les travaux antérieurs, posent une définition des espaces connectifs identique à celle de Börger, mais sans les morphismes. En 1998, Jean Serra élargit cette définition à celle de \og connections\fg\footnote{Notion qui, \textit{a priori}, n'a rien à voir avec les notions de connexion en géométrie différentielle.} sur un treillis \cite{Serra:1998, Serra:2000}, la définition classique correspondant au cas où le treillis considéré est celui des parties d'un ensemble. Nos propres travaux sur le sujet\footnote{Voir en particulier \cite{Dugowson:2003, Dugowson:2007c, Dugowson:201012, Dugowson:201203}} ont débuté en 2002, à l'occasion d'une réflexion sur la topologie du jeu de go. L'espace borroméen --- constitué de trois points globalement connectés sans qu'ils le soient deux à deux, de sorte que sa structure connective est précisément celle du n\oe ud borroméen, raison évidente de cette dénomination --- y est l'un des tous premiers exemples que nous considérons, remarquable pour ne pouvoir être défini ni à partir d'une structure topologique sur l'ensemble des trois points, ni à partir de celle d'un graphe à trois sommets.  Notons au passage que la notion de \textit{connective spaces}, introduite en 2006 par Joseph Muscat et David Buhagiar\cite{Muscat_Buhagiar:2006} de façon également indépendante des contributions antérieures, plus proche des espaces topologiques et de ce fait plus restreinte que celle de Börger, est de notre point de vue trop restrictive  car elle ne permet pas de rendre compte de la structure connective des entrelacs\footnote{Voir plus loin la remarque \ref{rmq Muscat Buhagiar}.}.

Grâce aux espaces de séparation que nous introduisons dans la  section \ref{section espaces de separation},  la représentation  par entrelacs  d'un espace connectif fini s'interprète comme un cas particulier de la notion de représentation d'un espace connectif dans un autre (section \ref{section representation connective}), généralisation que la consi\-dération des espaces connectifs infinis nous a conduit à développer afin de clarifier dans ce cas l'idée même de représentation.
Or, les premiers exemples de représentation d'espaces connectifs infinis (fibration de Hopf,  dynamiques de Lorenz, de Ghrist, etc.\footnote{Voir \cite{Dugowson:201203}, dont \cite{Dugowson:201112} constitue une version antérieure disponible en ligne.}) suggèrent fortement la nature dynamique de tels objets. L'idée de dynamique renvoyant, d'un point de vue topologique ou géométrique, à celle de feuilletage, cela nous a conduit à développer la notion de feuilletage connectif, objet de la section
\ref{section feuilletages connectifs}
du présent article. Les notions de feuilletages connectifs et de représentations connectives se révèlent alors, sous certaines conditions, adjointes l'une de l'autre. Une telle adjonction est traitée en section 
\ref{section adjonction representation et feuilletages}.
L'étude des feuilletages connectifs annonçant d'autres travaux portant plus spécifiquement sur l'aspect connectif des systèmes dynamiques\footnote{Voir \cite{Dugowson:201203} pour une première introduction aux dynamiques catégoriques connectives.}, il est devenu nécessaire d'élargir la notion d'ordre connectif, introduite dans \cite{Dugowson:201012} dans le cas des espaces finis, aux espaces connectifs infinis. C'est ce que nous faisons dans la section 
\ref{section ordre connectif},
où nous définissons en particulier l'ordre connectif d'un feuilletage connectif.

Enfin, un peu à part, la section \ref{section relations avec diffeologie}, qui aborde la question des relations fonctorielles entre espaces connectifs et espaces difféologiques, trouve néanmoins à s'articuler de plusieurs façons avec les thèmes précédents. En particulier, le fait que la théorie des espaces connectifs tende à s'orienter d'elle-même vers des notions dynamiques constitue une incitation à explorer 
les relations entre espaces connectifs et espaces difféo\-logiques, dans la mesure où, parce qu'elle généralise et unifie diverses  constructions liées à la géométrie différentielle, la difféologie devrait jouer un rôle croissant dans l'étude de diverses classes de systèmes dynamiques. Il pourrait en particulier être intéressant d'explorer le thème des feuilletages en adoptant un double point de vue connectif et difféologique. Nous ne le ferons pas ici, nous contentant d'une question plus élémentaire, celle de préciser, en mettant à jour les foncteurs en jeu dans ces questions, les conditions sous lesquelles un espace connectif est \og difféo\-logisable\fg, au sens où sa structure peut être associée à une structure difféologique sur le même ensemble de points, à l'exemple de la structure borroméenne, par laquelle nous avons vu qu'aura commencé l'exploration connective, et qui se rencontre aussi bien en difféologie.

\paragraph{Notations et rappels}\label{notations et rappels}

Conformément aux définitions et aux notations introduites dans notre article \cite{Dugowson:201012}, rappelons ou précisons notamment les points suivants :

\begin{itemize}
\item si $f$ est une application $A\rightarrow B$, nous notons $f$ ou $f^\mathcal{P}$ l'application de $\mathcal{P}A$ dans $\mathcal{P}B$ définie pour toute partie $U$ de $A$ par $f^\mathcal{P}(U)=\{f(u), u\in U\}$, 
\item si $\rho$ est une application $A\rightarrow \mathcal{P}B$, nous notons ${^\mu\rho}$ l'application de $\mathcal{P}A$ dans $\mathcal{P}B$ définie pour toute partie $U$ de $A$ par ${^\mu}\rho(U)=\bigcup_{u\in U} \rho(u)$,
\item pour toute catégorie $\mathbf{C}$, nous désignons par $\overrightarrow{\mathbf{C}}$ la classe de ses flèches, et par $\dot{\mathbf{C}}$ ou $\mathbf{C}_0$ la classe de ses objets. 
\item un espace connectif $X$ consiste en un couple noté $X=(\vert X\vert, \kappa(X))$, $\vert X\vert$ étant le support de l'espace, et $\kappa(X)$ sa structure,
\item on note $\mathbf{Cnc}$ la catégorie des espaces connectifs, et $\mathbf{Cnct}$ celle des espaces connectifs intègres, 
\item l'ensemble des structures connectives définies sur un ensemble constitue un treillis complet pour l'inclusion,
\item sur toute partie du support d'un espace connectif se trouve définie une structure connective dite \emph{induite}, qui est la moins fine faisant de l'injection canonique un morphisme connectif,
\item la structure connective engendrée par un ensemble de parties $\mathcal{A}$ d'un ensemble donné est notée $[\mathcal{A}]_0$, tandis que la structure connective intègre engendrée par $\mathcal{A}$ est notée  $[\mathcal{A}]_1$ ou simplement $[\mathcal{A}]$,
\item \label{systeme de generateurs} étant donné un espace connectif $X$, nous dirons en outre qu'un ensemble $\mathcal{A}$ de parties de $\vert X\vert$ constitue un système de générateurs de $X$ (ou de sa structure $\kappa(X)$) si $[\mathcal{A}]_0=\kappa(X)$,
\item $U_T:\mathbf{Top}\longrightarrow \mathbf{Cnct}$ désigne le \og foncteur d'oubli\fg\, de la catégorie des espaces topologiques dans la catégorie des espaces connectifs intègres.
\end{itemize}

\section{Espaces de séparation}\label{section espaces de separation}

Les espaces de séparation sont une manière de présenter les espaces connectifs intègres. Ils nous permettront d'interpréter les représentations par entrelacs comme des représentations connectives particulières (voir l'exemple \ref{representations par entrelacs comme representations}). Dans la suite, $E$ désigne un ensemble quelconque.

\begin{df}
On appelle \emph{dispositif de séparation} sur $E$ tout ensemble  $\mathcal{S}$ de paires $\{S,T\}$ de parties non vides et disjointes de $E$. 
 Les paires $\{S,T\}$ d'un tel dispositif sont appelés \emph{paires séparatrices}. 
\end{df}

\begin{df} Soit $\mathcal{S}$ un dispositif de séparation sur $E$. On dit qu'une partie $A$ de $E$ est \emph{séparée par $\mathcal{S}$}, et l'on note $(\mathcal{S}:A)$,
s'il existe dans $\mathcal{S}$ une paire séparatrice $\{S,T\}$ qui recouvre $A$ et dont chaque membre rencontre $A$ : 
$A\subset S\cup T$, $A\cap S\neq \emptyset$ et $A\cap T\neq \emptyset$.
\end{df}

\begin{rmq} Pour tout groupe $G$ de permutations de $E$ et tout dispositif de sé\-para\-tion $\mathcal{S}$ sur $E$, on a alors $(G\mathcal{S}:A)$ $\Leftrightarrow$ $\exists\varphi\in G$ tel que $(\mathcal{S}:\varphi(A))$, où $G\mathcal{S}$ est le dispositif de séparation défini par  $G\mathcal{S}=\{\{\varphi(S),\varphi(T)\}, \varphi\in G, (S,T)\in\mathcal{S}\}$. 
\end{rmq}

\begin{df} 
Soit $\mathcal{S}$ un dispositif de séparation sur $E$. L'ensemble $\kappa(\mathcal{S})=\{K\in \mathcal{P}(E), \neg (\mathcal{S}:K)\}$ 
constitue une structure connective intègre sur $E$. On note $E[\mathcal{S}]$  \emph{l'espace connectif défini sur $E$ par le dispositif de séparation
  $\mathcal{S}$}, de sorte que l'on a :
   $\vert E[\mathcal{S}] \vert=E$ et 
   $\kappa (E[\mathcal{S}])=\kappa(\mathcal{S})$.
\end{df}

\begin{thm} Tout espace connectif intègre peut être défini par un dispositif de séparation.
\end{thm}

{\noindent \textbf{Preuve}. On forme un dispositif de séparation adéquat en prenant tous les couples de parties disjointes non vides $(A,B)$ telles que toute \pagebreak[4] compo\-sante connexe de $A\cup B$ soit contenue soit dans $A$, soit dans $B$.\nolinebreak[4]\nopagebreak[4]
\begin{flushright}$\square$\end{flushright} }

\begin{exm} [Foncteur $V_T$] \label{exemple foncteur V_T}

On définit un foncteur $V_T:\mathbf{Top}\to\mathbf{Cnct}$ en posant, pour tout $X=(\vert X\vert,\mathcal{T}_X)\in \mathbf{Top}_0$ : $V_T(X)=\vert X\vert [\mathcal{S}_X]$, où $\mathcal{S}_X=\{\{S,T\}\in (\mathcal{T}\setminus\emptyset)^2, S\cap T=\emptyset\}$, et pour toute application continue $f:X\to Y$, $V_T(f)=f$.
En effet, une application continue transforme né\-ces\-saire\-ment une partie  non séparable par ouverts disjoints de l'espace de départ en partie non séparable par ouverts disjoints de l'espace d'arrivée.
Remarquons que toute partie connexe au sens topologique est néces\-sai\-re\-ment connexe au sens de cette nouvelle structure connective. Autrement dit, le foncteur $U_T$ est connectivement plus fin que le foncteur $V_T$. Par exemple, pour $\vert X\vert=\{1,2,3\}$ admettant pour ouverts non triviaux $\{1,2\}$ et $\{1,3\}$,  $\{2,3\}$ est non connexe dans $U_T(X)$ mais connexe dans $V_T(X)$. 
Par contre, dans un espace \emph{métrique} $X$, la connexité d'une partie est équivalente à l'impossibilité de la séparer par des ouverts \textit{disjoints} : $U_T(X)=V_T(X)$ (nous en donnons une preuve dans \cite{Dugowson:201203}, page 18).
\end{exm}

Une classe d'espaces de séparation importants est constituée des espaces affines réels : on munit l'espace affine réel  $E_n$ de dimension $n$ d'une structure connective notée $\sigma_n$ ou $\sigma$, et appelée la \emph{{structure connective usuelle de séparation}} sur $E_n$, en prenant pour dispositif de séparation $G\mathcal{S}$ avec $G$ le groupe des homéo\-mor\-phismes de l'espace topologique $E_n \simeq\mathbf{R}^n$ et pour $\mathcal{S}$ le singleton $\{\{S,T\}\}$, avec $S$ et $T$ les deux demi-espaces ouverts définis par un hyperplan quelconque.

\begin{df} L'espace connectif $(E_n,\sigma)$ est appelé \emph{l'espace usuel de séparation $n$-dimensionnel}.
\end{df}

\begin{prop}
La structure connective de l'espace usuel de séparation $(E_n,\sigma)$ n'est pas celle d'un espace topologique.
\end{prop}

\noindent \textit{Preuve}. On vérifie facilement que dans tout espace topologique, si $A$ et $B$ sont deux parties connexes non vides et que $x$ est un point de l'espace tel que $x \notin A\cup B$ et que $A\cup\{x\}$ et $B\cup\{x\}$ soient non connexes, alors $A\cup B \cup\{x\}$ est encore non connexe. Or, dans $(E_n,\sigma)$, si l'on prend par exemple pour $A$ une demi-sphère,  pour $B$ la demi-sphère complémentaire et pour $x$ le centre de la sphère $A\cup B$, la propriété précédente est contredite.
\begin{flushright}$\square$\end{flushright} 
\pagebreak[3]

\`{A} noter que l'espace connectif $(E_n,\sigma)$ est moins fin (il a plus de connexes) que l'espace connectif $(E_n,\tau)$ associé par $U_T$ à l'espace topologique usuel $E_n$.

On définit de même la sphère (réelle) usuelle de séparation $n$-dimen\-si\-onnelle
 $(S^n,\sigma)$  en prenant pour dispositif de séparation $G\mathcal{S}$ avec $G$ le groupe des homéo\-mor\-phismes de l'espace topologique $S^n$
  et pour $\mathcal{S}$ le singleton $\{\{S,T\}\}$, 
  avec $S$ et $T$ les deux ouverts séparés par une sphère $(n-1)$-dimensionnelle plongée dans $S^n$.

%
%
%
%

\begin{rmq}\label{rmq Muscat Buhagiar} La propriété des espaces topologiques utilisée dans la preuve ci-dessus pour montrer que les espaces $(E_n,\sigma)$ ne sont pas topologiques fait partie de celles incorporées par Muscat et Buhagiar dans leur défi\-nition des \textit{connective spaces} \cite{Muscat_Buhagiar:2006}. Pour cette raison, les espaces de Muscat et Buhagiar ne permettent pas de rendre compte de la structure connective des entrelacs (voir ci-dessous l'exemple \ref{representations par entrelacs comme representations}).
\end{rmq}

\section{Représentations connectives}\label{section representation connective}
 
Le théorème de Brunn-Debrunner-Kanenobu (voir \cite{Dugowson:201012, Dugowson:201203}) concerne la représentation \textit{par entrelacs} des espaces connectifs finis. Le point de vue que nous soutenons ici est que la représentation d'un espace connectif fini par entrelacs doit être comprise comme un cas particulier de la notion générale de représentation d'un espace connectif dans un autre, objet de la présente section.


L'idée des représentations connectives est d'associer à tout \textit{point} de l'espace représenté une \textit{partie} non vide de l'espace dans lequel a lieu la représentation. Pour cela, nous aurons besoin de faire appel au foncteur $\mathcal{P}^*$ ainsi défini :

\begin{df} On définit un endofoncteur $\mathcal{P}^*$ de la catégorie $\mathbf{Cnc}$, appelé \emph{puissance connective générale} ou \emph{espace connectif des parties non vides}, en associant à tout espace connectif $X$ l'espace connectif, noté $\mathcal{P}^*_X$ (ou $\mathcal{P}^*(X)$, ou $\mathcal{P}^*X$) défini par
\begin{itemize}
\item son support $\vert \mathcal{P}^*_X\vert =\mathcal{P}^*_{\vert X\vert}$,
\item et sa structure connective $\kappa(\mathcal{P}^*_X) =\{\mathcal{A}\in \mathcal{P}\mathcal{P}^*_{\vert X\vert}, \bigcup \mathcal{A} \in \kappa(X)\}$,
\end{itemize}
et en associant à tout morphisme connectif $f:X\rightarrow Y$, le morphisme connectif noté $f^\mathcal{P}$ ou simplement $f$, défini pour toute partie non vide $A$ de $\vert X\vert$ par $f(A)=\{f(a), a\in A\}$.  
\end{df}

\begin{rmq} L'espace connectif $\mathcal{P}^*X$ n'est pas intègre en général, même lorsque $X$ l'est. 
\end{rmq}

\begin{rmq} On définit de même un endofoncteur $\mathcal{K}^*$ de la catégorie $\mathbf{Cnct}$, appelé \emph{puissance connective intègre} ou \emph{espace des parties connexes non vides}, en associant à tout espace connectif intègre $X$ l'espace connectif $\mathcal{K}^*_X$ de ses parties connexes non vides, de  structure connective $\kappa(\mathcal{K}^*_X) =\{\mathcal{A}\in \mathcal{P}(\vert \mathcal{K}^*_X\vert), \bigcup \mathcal{A} \in \kappa(X)\}$, et à tout morphisme connectif $f:X\rightarrow Y$, le morphisme connectif encore noté $f^\mathcal{P}$ ou $f$, défini pour toute partie connexe non vide $K$ de $\vert X\vert$ par $f(K)=\{f(a), a\in K\}$.  
\end{rmq}




\begin{df} [Représentation connective] On appelle \emph{représentation} connective d'un espace connectif $X$ dans un espace connectif $Y$ tout morphisme connectif de $X$ dans l'espace connectif $\mathcal{P}^*(Y)$. On écrira  $\rho:X\rightsquigarrow Y$ pour exprimer que $\rho$ est une représentation de $X$ dans $Y$. Étant donnée $\rho$ une telle représentation, $Y$ sera appelé l'\emph{espace} de $\rho$, et sera noté $Y=sp(\rho)$; $X$ sera appelé l'\emph{objet} de $\rho$, et sera noté $ob(\rho)$.
\end{df}

Dans le cas où $X$ est intègre, une représentation $\rho:X\rightsquigarrow Y$ s'identifie à un morphisme connectif de $X$ dans l'espace intègre $\mathcal{K}^*(Y)$.

\begin{df} On dit qu'une représentation $f:X\rightsquigarrow Y$ est \emph{intègre} si son objet et son espace sont tous deux intègres.
\end{df}


Soit maintenant $\epsilon$ la transformation naturelle $Id_{\mathbf{Cnc}}\rightarrow \mathcal{P}^*$ définie pour tout espace connectif $X$ par $\forall x\in \vert X\vert, \epsilon_X(x)=\{x\}$, et $\mu$ la transformation naturelle
$\mathcal{Q}^*=\mathcal{P}^*\circ \mathcal{P}^*\rightarrow \mathcal{P}^*$ définie par $\forall \mathcal{A}\in \mathcal{Q}^*_{\vert X\vert}, \mu_X(\mathcal{A})=\bigcup\mathcal{A}$. Le triplet $(\mathcal{P}^*,\epsilon,\mu)$ constitue alors une monade sur $\mathbf{Cnc}$. La catégorie de Kleisli associée à cette monade a pour objets les espaces connectifs, et pour morphismes les re\-pré\-sen\-ta\-tions, la composée de deux représentations $\rho:X\rightsquigarrow Y$ et $\tau:Y\rightsquigarrow Z$ étant définie pour tout $x\in X$ par $\tau\odot \rho(x)=\mu_Z(\tau^\mathcal{P} (\rho(x)))\subset Z$.
%

Étant donnée une représentation $\rho:X\rightsquigarrow Y$, on notera 
${^\mu{\rho}}$ l'application de $\mathcal{P}^*X$ dans $\mathcal{P}^*Y$ définie par ${^\mu{\rho}}=\mu_Y\circ{\rho^\mathcal{P}}$. Une représentation de $X$ dans $Y$ est donc une application $\rho$ de $X$ dans $\mathcal{P}^*_Y$ telle que $^\mu \rho$ transforme toute partie connexe non vide de $X$ en une partie connexe non vide de $Y$.

Pour toute partie non vide $A$ de $X$, on a donc
\[ 
{^\mu{\rho}}(A)=\mu_Y(\rho^\mathcal{P}(A))=\mu_Y(\{\rho(a), a\in A\})=\bigcup_{a\in A}{\rho(a)}\subset Y
\]
tandis que la composée de deux représentations s'écrit
\[ 
\tau\odot \rho = {^\mu{\tau}}\circ \rho.
\]

\begin{rmq} En prenant  $\epsilon:Id_{\mathbf{Cnct}}\rightarrow \mathcal{K}^*$
 et $\mu:\mathcal{K}^*\mathcal{K}^*\rightarrow \mathcal{K}^*$ 
définis comme ci-dessus, le triplet $(\mathcal{K}^*,\epsilon,\mu)$ constitue de même une monade sur $\mathbf{Cnct}$, dont la catégorie de Kleisli associée a pour objets les espaces connectifs intègres, et pour morphismes les re\-pré\-sen\-ta\-tions intègres, avec la composition des re\-pré\-sen\-ta\-tions définie comme pour le cas général.
\end{rmq}

%






\begin{df}[Représentations claires et distinctes]
Soit $\rho:X\rightsquigarrow Y$ une représentation d'un espace $X$ dans un espace $Y$. On dit que $\rho$ est \emph{claire} si 
$\forall A\in \mathcal{P}_{\vert X\vert}, A\notin \kappa(X)\Rightarrow {^\mu\rho}(A)\notin \kappa(Y)$.
On dit que $\rho$ est \emph{distincte} si
$\forall(x,y)\in X^2, x\neq y\Rightarrow \rho(x)\cap \rho(y)=\emptyset$. 
\end{df}

%
%
%

\begin{exm} Une représentation claire et distincte de l'espace borroméen $\mathcal{B}_3$ est obtenue en associant à chacun de ses points une des trois composantes d'un \textit{noeud borroméen} plongé dans $(E_3,\sigma_3)$. Plus généralement, les entrelacs brunniens constituent, dans $(E_3,\sigma_3)$, des représentations claires et distinctes  des espaces connectifs brunniens. 
\end{exm}

On vérifie facilement la première partie du théorème suivant (voir \cite{Dugowson:201203}, théorème 10), et le théorème de Brunn-Debrunner-Kanenobu entraîne alors la seconde partie :

\begin{thm} \label{representation canonique} Tout espace connectif admet une représentation claire et distincte dans un espace intègre. En particulier, tout espace connectif fini admet une représentation par entrelacs, les points non connexes étant représentés par deux ou plusieurs composantes séparables de tels entrelacs.
\end{thm}

%
%
%


\begin{df} [Représentations de type $\mathcal{S}$]
Soit  $\mathcal{S}$ un dispositif de séparation sur un ensemble $Y$. On appelle \emph{représentation de type $\mathcal{S}$} toute représentation claire et distincte d'un espace connectif $X$ dans l'espace $Y[\mathcal{S}]$.
\end{df}

\begin{exm}\label{representations par entrelacs comme representations} Toute représentation par entrelacs d'un espace connectif est une représentation de type $\mathcal{S}$, où $\mathcal{S}$ est un dispositif de séparation engendrant l'espace de séparation usuel $(E_3,\sigma_3)$. 
\end{exm}



\begin{df}\label{categorie des representations et rcd}[Catégorie des représentations] On définit une catégorie $\mathbf{RC}$, dite \emph{catégorie des re\-pré\-sen\-ta\-tions connectives} en prenant pour objets les re\-pré\-sen\-ta\-tions connectives, et pour morphismes  d'une re\-pré\-sen\-ta\-tion $\rho:A\rightsquigarrow B$ vers une re\-pré\-sen\-ta\-tion $\rho':A'\rightsquigarrow B'$ les couples $(\alpha,\beta)$ où $\alpha:A\to A'$ et $\beta:B\to B'$ sont des morphismes connectifs tels que
$
\beta^{\mathcal{P}}\circ \rho \subset \rho'\circ\alpha,
$
 au sens où, pour tout $a\in A$, $\beta^{\mathcal{P}}(\rho(a))\subset \rho'(\alpha(a))$.
 La sous-catégorie pleine de $\mathbf{RC}$ admettant pour objets les représentations claires et distinctes sera notée $\mathbf{RCD}$.
\end{df}

\begin{exm}[points d'une représentation] La catégorie $\mathbf{RC}$ admet comme objet final l'unique représentation $\mathbf{1}_\mathbf{RC}:\bullet\mapsto\{\bullet\}$ d'un singleton connecté dans lui-même. Un point d'une représentation $\rho:A\rightsquigarrow B$ est alors un morphisme $\mathbf{1}_\mathbf{RC}\to\rho$, c'est-à-dire la donnée d'un point connecté $p$ de $A$ et d'un point connecté $q$ de $\rho(p)\subset B$. En particulier, si l'objet ou l'espace d'une représentation ne possède pas de point intègre, celle-ci n'a pas de point.
\end{exm}

\section{Feuilletages connectifs}\label{section feuilletages connectifs}


\begin{df} [Feuilletage connectif]
Un \emph{feuilletage connectif}  est un triplet $(E,\kappa_0, \kappa_1)$ constitué d'un ensemble $E$ appelé le \emph{support} du feuilletage, et d'un couple $(\kappa_0, \kappa_1)$ de structures connectives sur $E$, la première, $\kappa_0$, étant dite structure connective \emph{interne}, et la seconde, $\kappa_1$, structure connective \emph{externe}. Lorsque que $\kappa_0\subset\kappa_1$, le feuilletage est dit \emph{régulier}.
\end{df}

Lorsqu'une partie de $E$ est connexe pour $\kappa_0$ (resp. $\kappa_1$), on dit aussi qu'elle est $\kappa_0$-connexe, ou encore qu'elle est \textit{connexe interne}, ou encore \textit{intérieurement connexe} (resp. $\kappa_1$-connexe, ou \textit{connexe externe}, ou encore \textit{extérieurement connexe}). Étant donné un feuilletage connectif $Z$, on notera $\vert Z\vert$ son support, $\kappa_0(Z)$ sa structure connective interne et $\kappa_1(Z)$ sa structure connective externe, de sorte que $Z=(\vert Z\vert, \kappa_0(Z), \kappa_1(Z))$. Souvent, on notera  $Z_0$ l'espace connectif intérieur $Z_0=(\vert Z\vert, \kappa_0(Z))$, et $Z_1$ l'espace connectif extérieur $Z_1=(\vert Z\vert, \kappa_1(Z))$.

\begin{df}\label{categorie FC} La \emph{catégorie des feuilletages connectifs} $\mathbf{FC}$ a pour objets les feuilletages connectifs, et pour morphismes d'un feuilletage $Z$ vers un feuilletage $Z'$ les applications $\vert Z\vert\to\vert Z'\vert$ qui sont connectives de $Z_i=(\vert Z\vert,\kappa_i(Z))$ vers $Z'_i=(\vert Z'\vert,\kappa_i(Z'))$ pour chacun des deux indices $i\in\{0,1\}$.
\end{df}

%
%
%

\begin{df} [Feuilles] Soit $Z$ un feuilletage. On appelle \emph{domaine} de $Z$, et on note $dom(Z)$, la partie présente de la structure interne $\kappa_0(Z)$. On appelle \emph{feuilles} de $Z$ les composantes connexes non vides de la structure interne $\kappa_0(Z)$. La \emph{structure interne} d'une feuille $F$ est la structure connective induite sur $F$ par $\kappa_0(Z)$. La structure externe de $F$ est la structure induite sur $F$ par $\kappa_1(Z)$. 
\end{df}

Pour tout feuilletage $Z$, on note $\mathcal{F}(Z)$ l'ensemble des feuilles de $Z$. Si $dom(Z)$ est non vide, $\mathcal{F}(Z)$ en constitue une partition.  


\begin{rmq}Par définition, chaque feuille est intérieurement connexe. Par contre, si le feuilletage n'est pas régulier, une feuille peut ne pas être exté\-rieu\-rement connexe. 
\end{rmq}


\begin{df} On dira qu'un morphisme de feuilletages $\phi:Z\to Z'$ est \emph{strict} si $\phi^\mathcal{P}$ transforme toute feuille de $Z$ en une feuille de $Z'$. La catégorie ayant pour objets les feuilletages connectifs et pour morphismes les morphismes de feuilletages stricts sera notée $\mathbf{FS}$.
\end{df}


\begin{exm} Un espace topologique $Y$ muni d'une relation d'equivalence $\rho$ définit un feuilletage connectif, en prenant $(\vert Z\vert,\kappa_1(Z))=U_T(Y)$ et $\kappa_0(Z)=\rho$, la structure connective associée à la relation d'équivalence $\rho$. 
\end{exm}





\begin{exm}  \`{A} toute variété feuilletée on associe le feuilletage connectif régulier défini sur le même ensemble de points en prenant pour structure connective interne celle associée à la topologie la plus fine du feuilletage (celle de plus faible dimension), et pour structure connective externe celle associée à la topologie la moins fine du feuilletage.
\end{exm}


\begin{df}[Espace induit des feuilles]\label{espace induit des feuilles} Soit $Z=(\vert Z\vert,\kappa_0(Z), \kappa_1(Z))$ un feuilletage connectif. 
L'\emph{espace induit} des feuilles de $Z$ est l'espace connectif noté $\mathcal{F}^\downarrow(Z)$, ou plus simplement $Z^\downarrow$,  de support $\vert Z^\downarrow\vert$ l'ensemble $\mathcal{F}(Z)$ des feuilles de $Z$, 
et de structure connective celle qui y est induite par l'espace connectif des parties non vides $\mathcal{P}^*(Z_1)$, 
où $Z_1=(\vert Z\vert, \kappa_1(Z))$, de sorte qu'un ensemble $\mathcal{A}$ de feuilles 
est $\kappa(Z^\downarrow)$-connexe si et seulement si $\bigcup_{F\in\mathcal{A}}{F}\in \kappa_1(Z)$.
\end{df}

\begin{rmq} Si une $\kappa_0$-composante connexe n'est pas $\kappa_1$-connexe, elle définit un point non connexe de l'espace $Z^\downarrow$. Ainsi, l'espace induit des feuilles d'un feuilletage $Z$ est-il intègre si et seulement si toute composante connexe de la structure interne de $Z$ est extérieurement connexe.
\end{rmq}

\begin{rmq} 
Il existe une autre façon, qui pourrait d'ailleurs sembler plus naturelle, de munir l'espace des feuilles d'une structure connective. En effet, les feuilles d'un feuilletage étant
les composantes connexes de sa structure interne, elles sont également les classes d'une certaine relation d'équivalence partielle, d'où il découle très naturellement la définition de l'espace \textit{quotient} des feuilles. La structure quotient ainsi obtenue (appelée également structure \textit{sortante}) est moins fine que celle de l'espace induit (également appelé espace \textit{entrant}). Nous ne développerons pas ici la notion d'espace quotient des feuilles, renvoyant le lecteur intéressé aux sections §\,1.8 et §\,2.2.2. de  \cite{Dugowson:201203}.
\end{rmq}

\section{Une adjonction entre représentations et feuilletages}\label{section adjonction representation et feuilletages}

\subsection{Une famille de foncteurs $\mathbf{RC} \rightarrow \mathbf{FC}$}

\`{A} toute représentation connective $\rho:ob(\rho)\rightsquigarrow sp(\rho)$ on souhaite associer fonctoriellement un feuilletage $\Phi(\rho)$. Si, pour la structure externe du feuilletage, la structure de l'espace $sp(\rho)$ de la représentation s'impose, il y a par contre plusieurs choix possibles, \textit{a priori} légitimes,  pour la structure interne.  Pour préciser ces choix, nous aurons besoin de faire appel à ce que nous appellerons des structures connectives fonctorielles :

\begin{df} Une \emph{structure connective fonctorielle} est une application $\gamma$ définie sur la classe des espaces connectifs et qui à tout espace connectif $B$ associe une structure connective $\gamma(B)$ sur $\vert B\vert$ qui soit fonctorielle au sens où il existe un endofoncteur $\Gamma$ de $\mathbf{Cnc}$ défini 
sur les objets par $\Gamma(B)=(\vert B\vert,\gamma(B))$, 
et sur les flèches par $\Gamma(f)=f$.
Nous dirons qu'une structure connective fonctorielle $\gamma$ est plus fine qu'une autre, $\gamma'$, et l'on notera $\gamma\subset \gamma'$, si pour tout espace connectif $B$ on a $\gamma(B)\subset \gamma'(B)$.
\end{df}

Par exemple, notons respectivement $\kappa_D(B)$ et $\kappa_G(B)$ la structure connective désintégrée\footnote{C'est-à-dire la structure discrète non intègre, pour laquelle seul le vide est connexe.} et la structure connective grossière sur $\vert B\vert$.  Alors  $\kappa_D$ et  $\kappa_G$ sont des structures connectives fonctorielles. De même, l'application  $\kappa$ qui à tout espace connectif $B$ associe sa structure connective $\kappa(B)$ est une structure connective fonctorielle, et l'on a :
\[ 
\kappa_D\subset \kappa \subset \kappa_G.
 \]

Soit maintenant $(\gamma_0,\gamma_1)$ 
un couple de structures connectives fonctorielles tel que 
$\gamma_0 \subset \gamma_1$. \`{A} toute représentation connective 
$\rho : A \rightsquigarrow B$,
on associe le feuilletage $Z=\Phi_{(\gamma_0,\gamma_1)}(\rho)= \Phi(\rho)$ de support 
$\vert Z\vert=\vert B\vert$,
de structure externe $\kappa_1(Z)=\kappa(B)$ et de structure interne
\[ 
\kappa_0(Z)=[\bigcup_{i\in\{0,1\}}\bigcup_{a\in A_i}{( \gamma_i(B)\cap \mathcal{P}_{\rho(a)})}]_0,
 \]
où $A_0$ désigne la partie absente de $A$ et $A_1$ sa partie présente.

\begin{prop}\label{prop Phi foncteur} Soit $(\alpha,\beta):\rho\to\rho'$  un morphisme de représentations connectives. Alors l'application $\beta:\vert sp(\rho)\vert \to \vert sp(\rho')\vert$ est un morphisme de feuilletages $\beta:\Phi(\rho)\to\Phi(\rho')$.
\end{prop}
\noindent \textbf{Preuve}. Posons $A=ob(\rho)$, $B=sp(\rho)$, $A'=ob(\rho')$ et $B'=sp(\rho')$.
Si $K$ est une partie extérieurement connexe du feuilletage $Z=\Phi(\rho)$, alors $\beta^\mathcal{P}(K)$ est une partie extérieurement connexe de $Z'=\Phi(\rho')$ puisque les structures extérieures des feuilletages coïncident avec les structures des espaces de représentation que respecte $\beta$. 
Soit maintenant $K\in \bigcup_{i\in\{0,1\}}\bigcup_{a\in A_i}{( \gamma_i(B)\cap \mathcal{P}_{\rho(a)})}$. Il faut et il suffit de montrer que $\beta^\mathcal{P}(K)$ est intérieurement connexe dans $Z'$.
Si $K\in\gamma_0(B)\cap \mathcal{P}_{\rho(a)}$  avec $a\in A_0$, alors $\beta(K)\in\gamma_0(B')$ puisque $\gamma_0$ est fonctoriel. Et $K\subset \rho(a) \Longrightarrow \beta^\mathcal{P}(K)\subset \beta^\mathcal{P}(\rho(a))\subset \rho'(\alpha(a))$. Si $a'=\alpha(a)\in A'_0$, on a alors $\beta^\mathcal{P}(K)\in \gamma_0(B')\cap \mathcal{P}_{\rho(a')}\subset\kappa_0(Z')$, tandis que si $a'\in A'_1$, on a $\beta^\mathcal{P}(K)\in \gamma_1(B')\cap \mathcal{P}_{\rho(a')}\subset\kappa_0(Z')$,
puisque $\gamma_0\subset\gamma_1$.
Si $K\in\gamma_1(B)\cap \mathcal{P}_{\rho(a)}$ avec $a\in A_1$, on a nécessairement $a'=\alpha(a)\in A'_1$, et comme précédemment le fait que $\beta^\mathcal{P}(\rho(a))\subset \rho'(a')$ permet de conclure que $\beta^\mathcal{P}(K)\in \gamma_1(B')\cap \mathcal{P}_{\rho(a')}\subset\kappa_0(Z')$,
puisque $\gamma_1$ est fonctoriel.
\begin{flushright}$\square$\end{flushright} 
\pagebreak[3]

On en déduit que l'application $\Phi=\Phi_{(\gamma_0,\gamma_1)}$ qui à toute représentation $\rho$ associe le feuilletage $\Phi(\rho)$ et à tout morphisme de représentations $(\alpha,\beta)$ associe $\beta$ est un  foncteur $\mathbf{RC} \rightarrow \mathbf{FC}$. Dans le cas où $\gamma_0=\gamma_1=\gamma$, on le notera simplement $\Phi_\gamma$. Lorsque $\gamma=\kappa_G$ (resp. $\gamma=\kappa_D$), on notera simplement $\Phi_G$ (resp. $\Phi_D$) le foncteur $\Phi_\gamma$\label{notation PHIG}.

%


\begin{prop}\label{Phik regulier} Pour toute représentation connective $\rho$, le feuilletage $\Phi_\kappa (\rho)$ est régulier.
\end{prop}
\noindent \textbf{Preuve}. La structure interne de $\Phi_\kappa (\rho)$ étant, par définition, engendrée par des parties extérieurement connexes, elle est nécessairement plus fine que la structure externe.
\begin{flushright}$\square$\end{flushright} 
\pagebreak[3]

La proposition suivante découle immédiatement des définitions.

\begin{prop} Soient $(\gamma_0,\gamma_1)$ un couple de structures connectives fonctorielles, tel que $\gamma_0 \subset\gamma_1$, et soit $\Phi=\Phi_{(\gamma_0,\gamma_1)}$ le foncteur $\mathbf{RC}\to\mathbf{FC}$ associé. 
Si $\rho$ est une représentation distincte alors, pour $Z=\Phi(\rho)$, on a
\[ 
\kappa_0(Z)=\bigcup_{i\in\{0,1\}}\bigcup_{a\in A_i}{( \gamma_i(B)\cap \mathcal{P}_{\rho(a)})},
 \]
de sorte que si $\gamma_1\supset\kappa$ et que $a$ est un point connecté de $ob(\rho)$, alors $\rho(a)$ est 
une feuille de $Z$.
\end{prop}

\begin{cor}\label{cor feuilles representation} Si $\rho$ est une représentation distincte d'un objet intègre, les feuilles de $\Phi_\kappa (\rho)$ sont les parties de $sp(\rho)$ de la forme $\rho(a)$ :
\[ 
\mathcal{F}(\Phi_\kappa(\rho))=\{\rho(a), a\in ob(\rho)\}.
 \]
\end{cor}

\subsection{Le foncteur $\mathcal{R}^\downarrow:\mathbf{FC}\rightarrow \mathbf{RCD}$}

On définit de la manière suivante un foncteur $\mathbf{FC}\rightarrow \mathbf{RCD}$, appelé \emph{représentation induite} et noté $\mathcal{R}^\downarrow$.

\paragraph{Définition de $\mathcal{R}^\downarrow$  sur les objets.}
\`{A} tout feuilletage $Z$,  $\mathcal{R}^\downarrow$  associe la représentation $\mathcal{R}^\downarrow(Z):Z^\downarrow \rightsquigarrow Z_1$ de l'espace induit des feuilles $Z^\downarrow$, définie pour toute feuille $F\in \mathcal{F}(Z)=\vert Z^\downarrow \vert$ par
\[ 
\mathcal{R}^\downarrow(Z)(F)=F\subset\vert Z_1\vert.
 \]
$\mathcal{R}^\downarrow(Z)$ est bien une représentation connective puisque, par définition, un ensemble de feuilles est connexe dans $\mathcal{P}_{Z_1}$ si et seulement si son union est connexe dans $Z_1$, et que cette dernière propriété caractérise la structure connective de l'espace induit $Z^\downarrow$.  Il est en outre immédiat que la représentation $\mathcal{R}^\downarrow(Z)$ est claire (si un ensemble de feuilles est non connexe dans $Z^\downarrow$, alors leur union est également non connexe dans l'espace externe du feuilletage), et distincte (deux point différents, c'est-à-dire deux feuilles différentes, sont représentées par deux composantes connexes internes nécessairement disjointes).
\begin{prop}\label{objet integre si feuilletage regulier} Si le feuilletage $\mathcal{Z}$ est régulier, l'objet de la représentation $\mathcal{R}^\downarrow \mathcal{Z}$ est intègre.
\end{prop}
\noindent \textbf{Preuve}. Toute feuille étant extérieurement connexe, elle constitue un singleton connexe de $ob(\mathcal{R}^\downarrow \mathcal{Z})$.
\begin{flushright}$\square$\end{flushright} 
\paragraph{Définition de $\mathcal{R}^\downarrow$  sur les flèches.}
$\mathcal{R}^\downarrow$ est défini sur les flèches de $\mathbf{FC}$ en associant à tout morphisme de feuilletages $\phi:Z\to Z'$ le morphisme de représentations $(\phi_0,\phi_1)$, où $\phi_0: Z^\downarrow \to Z'^\downarrow$  est défini pour toute feuille $F\in \mathcal{F}(Z)$ par : $\phi_0(F)$ est celle des composantes connexes de l'espace interne $(\vert Z'\vert,\kappa_0(Z')$ qui contient le $\kappa_0(Z')$-connexe $\phi^\mathcal{P}(F)$, et où $\phi_1:Z_1\to Z'_1$ est le morphisme connectif qui en tant qu'application ensembliste coïncide avec $\phi$.

\begin{rmq}[Représentation quotient $\mathcal{R}^\uparrow$] \`{A} tout feuilletage on peut aussi associer fonctoriellement une représentation claire et distincte de son espace quotient de feuilles dans l'ensemble ambiant muni d'une structure connective adaptée. Le lecteur intéressé pourra se rapporter à la section §\,2.3.2. de \cite{Dugowson:201203}.
\end{rmq}

\subsection{L'adjonction $\mathcal{R}^\downarrow\dashv \Phi_\kappa $}

Pour établir cette adjonction entre les foncteurs $\mathcal{R}^\downarrow$ et $\Phi_\kappa$ lorsqu'ils sont restreints à certaines catégories de feuilletages et de représentations, nous faisons appel aux trois lemmes suivants (lemme \ref{lemme beta est un morphisme de feuilletages } à lemme \ref{lemme existence de alpha donnant representation}).

\begin{lm}\label{lemme beta est un morphisme de feuilletages } Soit $Z$ un feuilletage \emph{régulier}, $\rho$ une représentation quelconque, et $(\alpha,\beta) :\mathcal{R}^\downarrow Z \rightarrow \rho $ un morphisme de représentations. Alors $\beta$ est un morphisme de feuilletages $Z\rightarrow \Phi_\kappa (\rho)$.
\end{lm}
\noindent \textbf{Preuve}. Par définition d'un morphisme de représentations, $\beta$ est un morphisme connectif $sp(\mathcal{R}^\downarrow Z) \rightarrow sp(\rho)$, autrement dit un morphisme connectif $Z_1\rightarrow(\Phi_\kappa(\rho))_1$.
D'autre part, en appliquant le foncteur $\Phi_\kappa$ au morphisme $(\alpha,\beta)$ (proposition \ref{prop Phi foncteur}), on en déduit que $\beta$ est un morphisme de feuilletages $\Phi_\kappa(\mathcal{R}^\downarrow Z) \rightarrow \Phi_\kappa(\rho)$, donc en particulier un morphisme pour les structures internes
$(\Phi_\kappa(\mathcal{R}^\downarrow Z))_0 \rightarrow (\Phi_\kappa(\rho))_0$.
Mais $Z$ étant régulier, $\kappa_0(Z)\subset \kappa_0(\Phi_\kappa(\mathcal{R}^\downarrow Z))$. En effet, tout connexe intérieur est trivialement inclus dans une composante connexe intérieur et, par la régularité de $Z$, est aussi un connexe extérieur, de sorte que, par définition de la structure $\kappa_0(\Phi_\kappa(\mathcal{R}^\downarrow Z))$, se trouve bien appartenir à celle-ci.
Finalement, on a à la fois $\beta:Z_1\rightarrow(\Phi_\kappa(\rho))_1$ et $\beta:Z_0\rightarrow(\Phi_\kappa(\rho))_0$, autrement dit $\beta$ est bien un morphisme $Z\rightarrow \Phi_\kappa(\rho)$.
\begin{flushright}$\square$\end{flushright} 
\pagebreak[3]
\begin{lm}\label{lemme beta determine alpha} Soient $\mathcal{Z}$ un feuilletage connectif, $\rho$ une représentation connective \emph{distincte} et $(\alpha,\beta):\mathcal{R}^\downarrow(Z) \rightarrow \rho$ un morphisme de représen\-tations. Alors la connaissance de $\beta$ détermine celle de $\alpha$. Autrement dit, si $(\alpha',\beta):\mathcal{R}^\downarrow(Z) \rightarrow \rho$ est également un morphisme de représentations, on a nécessai\-rement $\alpha=\alpha'$.
\end{lm}
\noindent \textbf{Preuve}. 
Par définition, $\alpha$ est un morphisme connectif de 
$ob(\mathcal{R}^\downarrow Z)=\mathcal{F}^\downarrow(Z)=Z^\downarrow$ 
dans $ob(\rho)$. 
Soit $F\in ob(\mathcal{R}^\downarrow(Z))$, 
autrement dit une composante connexe de $Z_0=(\vert Z\vert, \kappa_0(Z))$. 
Par définition d'un morphisme de représentations, on a l'inclusion
$\beta^\mathcal{P}(\mathcal{R}^\downarrow Z (F))\subset \rho(\alpha(F))$.
Or, $\mathcal{R}^\downarrow Z (F)=F \subset \vert Z \vert$, d'où 
$\beta^\mathcal{P}(F)\subset\rho(\alpha(F))$. 
\pagebreak[3]
La représentation $\rho$ étant distincte, il n'y a au plus qu'un point $a$ de $ob(\rho)$ pouvant vérifier $\beta^\mathcal{P}(F)\subset\rho(a)$, d'où l'unicité annoncée. \begin{flushright}$\square$\end{flushright} 
\pagebreak[3]
\begin{lm}\label{lemme existence de alpha donnant representation}
Soit $Z$ un feuilletage, $\rho$ une représentation claire et distincte, d'objet $ob(\rho)$ intègre, et soit $\beta:Z\rightarrow \Phi_\kappa(\rho)$ un morphisme de feuilletages. Alors il existe un et un seul morphisme connectif $\alpha:\mathcal{F}^\downarrow Z \rightarrow ob(\rho)$ tel que $(\alpha,\beta)$ soit un morphisme de représentations $\mathcal{R}^\downarrow Z\rightarrow\rho$.
\end{lm}
\noindent \textbf{Preuve}. S'il existe, le morphisme $\alpha$ est unique d'après le lemme \ref{lemme beta determine alpha}. Précisons l'application ensembliste $\alpha:\mathcal{F}Z\rightarrow\vert ob(\rho)\vert$  dont, néces\-saire\-ment, il s'agit. Pour $F\in \mathcal{F}Z$, on a $\beta^\mathcal{P}(F)\in \kappa_0(\Phi_\kappa(\rho))$, puisque $\beta$ préserve aussi les morphismes internes. 
Notons $\overline{\beta^\mathcal{P}(F)}$ la composante $\kappa_0(\Phi_\kappa(\rho)$-connexe contenant $\beta^\mathcal{P}(F)$. Alors $\overline{\beta^\mathcal{P}(F)}\in\mathcal{F}(\Phi_\kappa(\rho))$. D'après le corollaire\,\ref{cor feuilles representation}, il existe alors un élément unique $a_F\in ob(\rho)$ tel que $\overline{\beta^\mathcal{P}(F)}=\rho(a_F)$. L'application $\alpha$ est donc définie par $\alpha(F)=a_F$. Autrement dit, 
\[ 
\alpha(F)=a 
\Leftrightarrow \beta^\mathcal{P}(F) \subset \rho(A) 
\Leftrightarrow 
\beta^\mathcal{P}(F)\subset \overline{\beta^\mathcal{P}(F)}= \rho(A).
 \]
Il s'agit de prouver que l'application $\alpha$ ainsi définie est un morphisme connectif $\mathcal{F}^\downarrow Z \rightarrow ob(\rho)$, et que le couple $(\alpha,\beta)$ est bien un morphisme de repré\-senta\-tions. 
Soit donc $\mathcal{L}$ un ensemble $\kappa(Z^\downarrow)$-connexe de feuilles. Par définition de $Z^\downarrow$, on a $\bigcup_{F\in\mathcal{L}} F \in \kappa_1(Z)$, donc l'ensemble $W=\bigcup_{F\in\mathcal{L}} \beta^\mathcal{P}(F)$ vérifie $W \in \kappa_1(\Phi_\kappa(\rho))$.
Posons  \[\mathcal{A}=\alpha^\mathcal{P}(\mathcal{L})=\{a\in ob(\rho), \exists F\in\mathcal{L}, \rho(a) \supset \beta^\mathcal{P}(F)\}. \]
On veut montrer que $\mathcal{A}$ est une partie connexe de $ob(\rho)$. 
Or, $\rho$ étant claire, il suffit pour cela de prouver que $^\mu \rho (\mathcal{A})= 
\bigcup_{F\in\mathcal{L}}\overline{\beta^\mathcal{P}(F)}$ est connexe dans $sp(\rho)$. 
Par définition, les $\overline{\beta^\mathcal{P}(F)}$ sont $\kappa_0(\Phi_\kappa(\rho))$-connexes. Mais, le feuilletage $\Phi_\kappa(\rho)$ étant régulier (proposition \ref{Phik regulier}), les $\overline{\beta^\mathcal{P}(F)}$ sont également $\kappa_1(\Phi_\kappa(\rho))$-connexes. Il en découle que
\[ 
\bigcup_{F\in\mathcal{L}}\overline{\beta^\mathcal{P}(F)}
=
\bigcup_{F\in\mathcal{L}}(\overline{\beta^\mathcal{P}(F)}\cup W)
 \]
\noindent est l'union de $\kappa_1(\Phi_\kappa(\rho))$-connexes d'intersection non vide. 

Ainsi, $\bigcup_{F\in\mathcal{L}}\overline{\beta^\mathcal{P}(F)}$ est $\kappa_1(\Phi_\kappa(\rho))$-connexe, autrement dit $\kappa_1(sp(\rho))$-connexe, de sorte que $\mathcal{A}$ est connexe dans $ob(\rho)$. 
Reste à vérifier que $\beta^\mathcal{P} \circ \mathcal{R}^\downarrow Z \subset \rho \circ \alpha$, mais c'est là une conséquence immédiate de la construction même de $\alpha$.
\begin{flushright}$\square$\end{flushright} 
\pagebreak[3]

Soit maintenant $\mathbf{FR}$ la sous-catégorie pleine de  $\mathbf{FC}$  constituée des feuilletages connectifs réguliers, et soit $\mathbf{RIO}$ la sous-catégorie pleine de $\mathbf{RCD}$ consti\-tuée des représentations claires et distinctes dont l'objet est intègre. Reprenons les notations  $\mathcal{R}^\downarrow$ et $\Phi_\kappa$ employées précédemment, mais pour désigner cette fois les restrictions de ces foncteurs à $\mathbf{FR}$ et à $\mathbf{RIO}$. 
D'après la proposition \ref{objet integre si feuilletage regulier}, on obtient bien de cette manière un foncteur $\mathcal{R}^\downarrow:\mathbf{FR}\rightarrow \mathbf{RIO}$. Et d'après la proposition \ref{Phik regulier}, on obtient de même un foncteur $\Phi_\kappa:\mathbf{RIO}\rightarrow \mathbf{FR}$.

Soit  $Z$ un feuilletage régulier, et $\rho$ une représentation claire et distincte d'un objet intègre. \`{A} tout morphisme de représenta\-tions $(\alpha,\beta):\mathcal{R}^\downarrow Z \rightarrow \rho$, on associe, d'après le lemme \ref{lemme beta est un morphisme de feuilletages }, le morphisme de feuilletages $\beta: Z\rightarrow \Phi_\kappa (\rho)$. Réciproquement, à tout morphisme de feuilletages $\beta: Z\rightarrow \Phi_\kappa (\rho)$, on associe d'après le lemme \ref{lemme existence de alpha donnant representation}, un unique morphisme de représentations $(\alpha,\beta):\mathcal{R}^\downarrow Z \rightarrow \rho$. On a ainsi construit des applications réciproques, donc bijectives, entre $Hom_{\mathbf{RIO}}(\mathcal{R}^\downarrow Z, \rho)$ et $Hom_{\mathbf{FR}}(Z, \Phi_\kappa(\rho))$, et le lecteur pourra vérifier que ces bijections sont naturelles par rapport à $Z$ et $\rho$. On peut ainsi énoncer :

\begin{thm}\label{Adjonction entre feuilletages et representations} Le  foncteur $\mathcal{R}^\downarrow:\mathbf{FR}\rightarrow \mathbf{RIO}$ est adjoint à gauche du foncteur $\Phi_\kappa:\mathbf{RIO}\rightarrow \mathbf{FR}$ :
\[ \mathcal{R}^\downarrow\dashv \Phi_\kappa \] 
\end{thm}

\begin{rmq} Par composition, les divers foncteurs considérés plus haut entre catégories de feuilletages connectifs et catégories de représentations connectives donnent lieu à d'autres foncteurs intéressants. Par exemple, notant $\rho^{\downarrow}_{G} = \mathcal{R}^\downarrow(\Phi_{G}(\rho) )$ 
la représentation associée à une représentation $\rho$ par l'endofoncteur $\mathcal{R}^\downarrow \circ\Phi_{G}$, on démontre\footnote{Voir \cite{Dugowson:201203}, proposition 18.} la proposition suivante :
\begin{prop} Si $\rho$ est une représentation claire et distincte, alors le couple d'applications $(\alpha,\beta)$ défini par 
$\alpha(a)=\rho(a)\in ob(\rho^{\downarrow}_{G})$ et $\beta=Id_{sp(\rho)}$ constitue un isomorphisme entre les représentations  $\rho$ et
$\rho^{\downarrow}_{G}$.
\end{prop}
%
\end{rmq}

\section{Ordre d'un espace connectif}\label{section ordre connectif}

 
On note $Ord$ la classe des ordinaux, $\omega_0$ ou $\aleph_0$ le plus petit ordinal infini, et $\aleph_1$ le plus petit ordinal non dénombrable, \textit{i.e.} l'ensemble des ordinaux dénombrables. Pour tout ordinal $\alpha$, nous notons en outre $\alpha^-$ l'ordinal défini par $\alpha^-=\beta$ si $\beta$ est prédécesseur de $\alpha$, et $\alpha^-=\alpha$ si $\alpha$ n'a pas de prédécesseur.

\begin{df} Soit $\alpha\in Ord$ un ordinal. Un ensemble (partiellement) ordonné $(R,\preceq)$ est dit \emph{supérieur ou égal} à $\alpha$, et l'on note $\alpha\leq R$, s'il existe une application strictement croissante de 
$\alpha$ dans 
$(R,\preceq)$.
\end{df}
Bien entendu, la définition précédente est compatible avec la relation d'ordre entre ordinaux. Soit maintenant $(R,\preceq)$ un ensemble ordonné. La classe des ordinaux $\alpha$ tels que $\alpha\leq R$ est bornée (en fonction du cardinal de $R$), c'est donc un ensemble, et c'est un ordinal puisque $\alpha\leq R\Rightarrow \beta\leq R$ pour tout $\beta\leq\alpha$. 

\begin{df} On appelle \emph{hauteur} de l'ensemble partiellement ordonné ${(R,\preceq)}$, et on note $\Gamma(R)$, l'ordinal
\[\Gamma(R)=\{\alpha\in Ord, \alpha\leq R\} \]
\end{df}

\begin{exm} $\mathbf{R}$ désignant
la droite réelle munie de l'ordre usuel, on a\footnote{On trouvera une preuve de ce fait dans \cite{Dugowson:201203}.} $\Gamma(\mathbf{R})=\aleph_1$.  
\end{exm}


Rappelons qu'une partie connexe non vide $K\in\kappa_X$ d'un espace connectif $X=(\vert X\vert,\kappa_X)$ est dite \emph{irréductible} si et seulement si elle n'appartient pas à la structure connective engendrée par les autres parties (voir \cite{Dugowson:201012}, section 2.2). La définition du graphe générique d'un espace connectif, donnée dans le cas fini dans \cite{Dugowson:201012} (section 7.1), s'étend immédiatement à tout espace connectif :

\begin{df}[Graphe générique]  Soit $X=(\vert X\vert,\kappa_X)$  un espace connectif.  On appelle graphe générique de $X$, et l'on note $(G_X,\subset)$, l'ensemble ordonné par l'inclusion des parties connexes irréductibles de $X$.
\end{df}

%

%
%


\begin{df}
Soit $X$ un espace connectif. On appelle \emph{ordre connectif} de $X$ l'ordinal $\Omega(X)=\Gamma(G_X)^{--}=\{\alpha\in Ord, \alpha+2\leq G_X\}$.  
\end{df}
%

Bien entendu, comme on le vérifie facilement, l'ordre connectif $\Omega(X)$ d'un espace connectif fini intègre $X$ coïncide avec l'ordre connectif défini dans \cite{Dugowson:201012}. Plusieurs exemples d'ordres connectifs infinis sont donnés dans \cite{Dugowson:201203}.

\`{A} tout entrelacs, qu'il comporte ou non un nombre fini de composantes, se trouve associé un espace connectif (voir l'exemple 3 dans \cite{Dugowson:201012} ou l'exemple 4 dans \cite{Dugowson:201203}). La notion d'ordre connectif conduit dès lors à la définition d'un nouvel invariant d'entrelacs  :
\begin{df} \emph{L'ordre connectif d'un entrelacs} est l'ordre connectif de l'espace connectif associé à cet entrelacs.
\end{df}

\begin{exm} L'ensemble des (classes d'équivalence d') entrelacs finis régu\-liers dans $\mathbf{R}^3$ étant dénombrable, on construit facilement un entrelacs dans $\mathbf{R}^3$ qui réalise l'union disjointe de tous les entrelacs finis réguliers. L'ordre connectif de l'entrelacs obtenu est $\omega_0$.
\end{exm}



La définition suivante est appelée à jouer un rôle important dans l'étude des dynamiques catégoriques connectives, puisqu'elle permet de définir l'ordre connectif d'une telle dynamique (voir \cite{Dugowson:201203}) :

\begin{df}[Ordre d'un feuilletage connectif] On appelle ordre, ou ordre connectif, d'un feuilletage connectif $Z$ l'ordre connectif de son espace induit de feuilles $Z^\downarrow$.
\end{df}

\section{Relations avec les espaces difféologiques}\label{section relations avec diffeologie}

On se propose dans cette section de préciser certaines relations fonctorielles entre espaces connectifs et espaces difféologiques, et en particulier de caractériser les espaces connectifs difféolo\-gisables. Nous commençons, après quelques rappels terminologiques, par préci\-ser la notion de difféologisabilité d'un espace connectif en définissant un foncteur d'oubli $U_{DC}$ de la catégorie des espaces difféologiques dans celle des espaces connectifs, puis nous donnons des conditions nécessaires de difféologisabilité d'un espace connectif avant de montrer, grâce à certains foncteurs de difféologisation, que ces conditions sont en fait suffisantes. Nous montrons ensuite que l'un de ces foncteurs est adjoint à droite du foncteur d'oubli, avant de conclure avec quelques remarques sur la notion d'application localement connective.


Pour tout ce qui concerne la difféologie, nous renvoyons à l'ouvrage \cite{PIZ:2013} de Patrick Iglesias-Zemmour. Rappelons néanmoins ici quelques notions et notations :

\begin{itemize}
\item on note $Param(E)$  l'ensemble des paramétrisations d'un ensemble $E$
(\cite{PIZ:2013}, art. 1.3)
;
\item une difféologie sur $E$ est une partie $\mathcal{D}\subset Param(E)$ vérifiant certains axiomes; les éléments de $\mathcal{D}$ s'appellent les plaques (\emph{plots}), une plaque $p$ de l'espace difféologique $(E,\mathcal{D})$ s'identifiant à une application lisse (\emph{smooth map}) définie sur un ouvert $U_p$ d'un espace de la forme $\mathbf{R}^{n_p}$ et à valeur  dans $E$; en particulier, une plaque définie sur $\mathbf{R}$ s'appelle un chemin (ou un chemin lisse), et l'ensemble des chemins dans $(E,\mathcal{D})$ est noté $Paths(E,\mathcal{D})$;
\item pour tout ensemble $\mathcal{L}$ de paramétrisations  $p:\mathbf{R}^{n_p}\supset U_p \rightarrow E$, autrement dit d'applications $p$ dans $E$, chacune étant définie sur un ouvert $U_p$ d'un espace de la forme $\mathbf{R}^{n_p}$ avec $n_p$ un entier naturel, on note $<\mathcal{L}>$ la difféologie engendrée par $\mathcal{L}$, c'est-à-dire la difféologie la plus fine sur $E$ contenant $\mathcal{L}$ (\cite{PIZ:2013}, art. 1.66);
\item on dit d'un ensemble de paramétrisations $\mathcal{L}$ de $E$ qu'il couvre $E$ si pour tout $a\in E$, il existe $p\in\mathcal{L}$ telle que $a\in p(U_p)=val(p)$;
\item une application $\sigma:\mathbf{R}\rightarrow E$ est dite \emph{stationnaire aux bords} s'il existe $\epsilon>0$ tel que pour tout $t\in]-\infty,\epsilon[$, $\sigma(t)=\sigma(0)$, et pour tout $t\in]1-\epsilon,+\infty[$, $\sigma(t)=\sigma(1)$; on note $stPaths(E,\mathcal{D})$ l'ensemble des chemins $\sigma\in Paths(E,\mathcal{D})$ qui sont stationnaires aux bords (voir \cite{PIZ:2013}, art. 5.4). 
\item une partie $A$ d'un espace difféologique $(E,\mathcal{D})$ est dite connectée si pour tout couple $(a_0,a_1)$ de points de $A$, il existe un chemin $\sigma\in Paths(E,\mathcal{D})$ reliant $a_0$ à $a_1$ en restant dans $A$ au sens où : $\sigma(0)=a_0$, $\sigma(1)=a_1$ et, pour tout $t\in \mathbf{R}$, $\sigma(t)\in A$ (voir \cite{PIZ:2013}, art. 5.9); par \og smashisation\fg\, (\cite{PIZ:2013}, art. 5.5.), on peut remplacer $Paths(E,\mathcal{D})$ par $stPaths(E,\mathcal{D})$ dans la définition des parties connectées d'un espace difféologique.
\end{itemize}


\begin{prop} L'ensemble $\mathcal{K}_\mathcal{D}$ des \emph{parties connectées} 
d'un espace difféologique
 $(E,\mathcal{D})$ 
constitue une structure connective intègre sur l'ensemble $E$.
\end{prop}

\noindent \textbf{Preuve}. Étant donnée $(K_i)_{i\in I}$ une famille de parties connectées de $(E,\mathcal{D})$, telle que $\bigcap_{i\in I} K_i\neq \emptyset$, deux points quelconques $a_1$ et $a_2$ de $L=\bigcup_{i\in I} K_i$ peuvent toujours être reliés par un chemin $\sigma\in Paths(E,\mathcal{D})$ : en effet, il existe un élément $a_0\in  \bigcap_{i\in I} K_i$ et, pour chaque $k\in\{1,2\}$, un chemin $\sigma_k\in\mathcal{D}$ reliant $a_k$ à $a_0$ en restant dans $L$. La \emph{smashed concatenation}\footnote{Voir \cite{PIZ:2013}, art. 5.5.} de $\sigma_1$ et $\sigma_2$ produit alors le chemin lisse $\sigma$ annoncé. En outre, il est clair que tout singleton est une partie connectée, puisque par définition d'une difféologie toutes les paramétrisations constantes, en particulier celles définies sur $\mathbf{R}$, appartiennent à $\mathcal{D}$.
\begin{flushright}$\square$\end{flushright} 

La proposition précédente, puisque par ailleurs toute application lisse entre espaces difféologiques transforme les parties connectées du premier en parties connectées du second (\cite{PIZ:2013}, art. 5.9), permet de définir un foncteur d'oubli\footnote{Au sens où il s'agit d'un foncteur fidèle entre catégories concrètes.} $U_{DC}$ de la catégorie $\mathbf{Diff}$  des espaces difféo\-lo\-giques  dans celle, $\mathbf{Cnct}$, des espaces connectifs intègres, en posant :
\begin{itemize}
\item pour tout espace $(E,\mathcal{D})\in \mathbf{Diff}_0$ : $U_{DC}(E,\mathcal{D})=(E,\mathcal{K}_\mathcal{D})$,
\item pour tout morphisme $f\in \overrightarrow{\mathbf{Diff}}$ :  $U_{DC}(f)=f$.
\end{itemize}

Dans la suite, nous étendons l'usage de l'expression  $U_{DC}$, permettant que la structure connective $\mathcal{K}_\mathcal{D}$ associée à la difféologie $\mathcal{D}$ soit également notée $U_{DC}(\mathcal{D})$. Ainsi, avec cette convention, on a, pour tout espace difféologique $(E,\mathcal{D})$,
\[U_{DC}(E,\mathcal{D})=(E,U_{DC}(\mathcal{D})).\]

\begin{df} Un espace connectif $(E,\mathcal{K})$ est dit \emph{difféologisable} s'il existe une structure difféologique $\mathcal{D}$ sur $E$ telle que $U_{DC}(\mathcal{D})=\mathcal{K}$. On notera $\mathbf{Cncd}$ la sous-catégorie pleine de $\mathbf{Cnct}$ ayant pour objets les espaces connectifs difféologisables.
\end{df}


\begin{lm}[Engendrement des connexes d'un espace difféologique]\label{lemme engendrement des connectes de diffeo par les chemins} Pour tout espace difféologique $(E,\mathcal{D})$, on a 
\[U_{DC}(\mathcal{D})=[\{\sigma([0,1]), \sigma\in stPaths(E,\mathcal{D})\}].\]
\end{lm}
\noindent \textbf{Preuve}. 
 Les chemins lisses étant des plaques, la connexité de l'intervalle réel $[0,1]$ entraîne celle des $\sigma([0,1])$, d'où $\mathcal{G}\subset \mathcal{K}_\mathcal{D}$, ce qui implique $[\mathcal{G}]\subset \mathcal{K}_\mathcal{D}$.  
 Soit maintenant une partie connectée non vide quelconque $K\in \mathcal{K}_\mathcal{D}$ et $a_0\in K$. Par définition de $\mathcal{K}_\mathcal{D}$, il existe, pour tout   $a\in K$,
un chemin  $\sigma_a$ 
tel que $\sigma_a(0)=a_0$, $\sigma_a(1)=a$ et 
$\sigma_a(\mathbf{R})\subset K$, chemin que, par \og \emph{smashisation} \fg (\cite{PIZ:2013}, art. 5.5) on peut prendre stationnaire aux bords : $\sigma_a\in stPaths(E,\mathcal{D})$. On a alors $K=\bigcup_{a\in K}{\sigma_a([0,1])}$, mais puisque $\bigcap_{a\in K}{\sigma_a([0,1])}\supset\{a_0\}\neq\emptyset$, cela prouve que \[K\in [\{\sigma([0,1]), \sigma\in stPaths(E,\mathcal{D})\}].\] 
\begin{flushright}$\square$\end{flushright}

\begin{lm} [Chemins d'une difféologie engendrée]\label{lemme caracterisation chemins en diffeologie engendree} Soit $\mathcal{L}$ un ensemble couvrant de paramétrisations de $E$, et soit $\sigma:\mathbf{R}\rightarrow E$ une paramétri\-sation définie sur $\mathbf{R}$. On a  $\sigma\in {stPaths(E,<\mathcal{L}>})$  si et seulement si les deux conditions suivantes sont satisfaites :
\begin{enumerate}
\item $\sigma$ est stationnaire au bord, 
\item il existe un entier $n\geq 1$ et une suite finie $(]a_k,b_k[, p_k,q_k)_{k\in \{1,...,n\}}$, avec  $a_k$ et $b_k$ des réels vérifiant 
\[(a_1<0)\,\,\mathrm{et}\,\, (\forall k\in\{1,..., n-1\}, a_k<a_{k+1}<b_k<b_{k+1})\,\,\mathrm{et}\,\, (b_n>1),\]
et telle que pour tout $k\in\{1,..., n\}$, on ait
\begin{itemize}
\item $(p_k:\mathbf{R}^{n_k}\supset U_k\rightarrow E)\in\mathcal{L}$, 
\item $q_k \in C^\infty(]a_k,b_k[,U_k)$, 
\item $\sigma_{\vert ]a_k,b_k[ }=p_k\circ q_k$,
\end{itemize}
où $n_k={n_{p_k}}$ désigne la dimension de la paramétrisation $p_k$, $U_k=U_{p_k}$ est le domaine de $p_k$ et $\sigma_{\vert ]a_k,b_k[ }$ désigne la restriction de $\sigma$ à $]a_k,b_k[$.
\end{enumerate}
\end{lm}
\noindent \textbf{Preuve}. D'après les axiomes qui définissent une difféologie, les conditions données sont clairement suffisantes pour avoir $\sigma\in <\mathcal{L}>$, les $p_k\circ q_k$ étant lisses par composition, et $\sigma$ étant lisse, puisque localement lisse, sur $[0,1]$ et constante, donc lisse, sur $]-\infty,\epsilon[$ et sur $]1-\epsilon,+\infty[$ pour un certain $\epsilon>0$. Stationnaire aux bord, $\sigma$ est donc bien dans $stPaths(E,<\mathcal{L}>)$. Réciproquement, étant donné un chemin stationnaire aux bords $\sigma\in stPaths(E,<\mathcal{L}>)$, la caractérisation des plaques d'une difféologie engendrée par une famille couvrante de paramétrisations donnée en \cite{PIZ:2013} (art. 1.68) implique que, pour tout $t\in \mathbf{R}$, il existe un voisinage ouvert $V_t\subset\mathbf{R}$ de $t$, une application de classe $C^\infty$ $q_t:V_t\rightarrow \mathbf{R}$ et une paramétrisation $p_t:\mathbf{R}\rightarrow E$ appartenant à $\mathcal{L}$ tels que $\sigma_{\vert V_t}=p_t\circ q_t$. Par restriction, on peut remplacer les ouverts $V_t$ par des intervalles ouverts $J_t\ni t$. La famille $(J_t)_{t\in [0,1]}$ est alors un recouvrement ouvert du compact $[0,1]$, on peut donc en extraire un sous-recouvrement fini, qu'après ré-indexation nous notons $(J_m)_{m\in\{1,..., N\}}$. Par récurrence finie, on construit alors de la façon suivante la suite des intervalles $]a_k,b_k[=I_k$ annoncés : parmi les intervalles $J_m$ qui contiennent $0$, on prend celui dont la borne supérieure est maximale, cela nous donne $I_1$, et l'on continue ainsi : ayant choisi les intervalles $I_1=]a_1,b_1[,...,I_k=]a_k,b_k[$, si $b_k>1$, on pose $n=k$ et l'on s'arrête, sinon on considère, parmi les intervalles $J_m$ qui contiennent $b_k$, celui dont la borne supérieure est maximale, ce qui nous donne $I_{k+1}$. Cette construction se poursuit tant que $b_k\leq 1$, puisqu'il existe alors un intervalle $J_m$ contenant $b_k$, mais elle s'achève nécessairement en un nombre fini $n$ d'étapes, d'où l'existence de $n$ tel que $b_n>1$. On vérifie alors aisément que les réels $a_k$ et $b_k$ ainsi obtenus satisfont les inégalités indiquées, d'où le résultat.
\begin{flushright}$\square$\end{flushright}

\begin{prop}\label{prop conditions necessaires de diffeologisabilite} Pour tout espace difféologique $(E,\mathcal{D})$, la structure connective intègre $U_{DC}(\mathcal{D})=\mathcal{K}_\mathcal{D}$ admet un système de générateurs
\footnote{Rappelons qu'un système de générateurs $\mathcal{G}$ d'une structure connective $\mathcal{K}$ est une partie $\mathcal{G}\subset\mathcal{K}$ telle que la structure connective engendrée par $\mathcal{G}$ vérifie $[\mathcal{G}]_0=\mathcal{K}$.}
 $\mathcal{G}\subset\mathcal{K}_\mathcal{D}$ tel que
\[\forall G\in\mathcal{G}, card(G)\leq \mathfrak{c},\]où $\mathfrak{c}$ désigne la puissance du continu.
\end{prop}

\noindent \textbf{Preuve}. Posons $\mathcal{G}=\{\sigma([0,1]), \sigma\in stPaths(E,\mathcal{D})\}$. Le lemme \ref{lemme engendrement des connectes de diffeo par les chemins} \pagebreak[3] dit que $\mathcal{G}$ est un système de générateurs de $\mathcal{K}_\mathcal{D}$, et puisque $card([0,1])=\mathfrak{c}$, 
tout  élément $G\in\mathcal{G}$ vérifie la condition de cardinalité indiquée.
\nopagebreak[3]
\begin{flushright}$\square$\end{flushright}

Pour toute paramétrisation $(p:R^{n_p}\supset U_p \rightarrow E)\in Param(E)$, nous noterons $\mathcal{KA}_p$ la structure connective sur $U_p$ constituée des parties de $U_p$ connexes par arcs, 
$\mathcal{T}_p$ la topologie sur $U_p$ induite par la topologie usuelle de $R^{n_p}$, et $\mathcal{KT}_p=U_{T}(\mathcal{T}_p)$ la structure connective sur $U_p$ associée à $\mathcal{T}_p$ par le foncteur d'oubli $U_{T}:\mathbf{Top}\rightarrow\mathbf{Cnct}$. On a donc, pour toute paramétrisation $p$, $\mathcal{KA}_p \subset \mathcal{KT}_p$.

\begin{prop}[Trois procédés de difféologisation des espaces connectifs\footnote{Une partie de ce qui est avancé ici nous a été suggéré par Anatole Khélif.}]\label{prop trois diffeologisations} Soit $(E,\mathcal{K})$ un espace connectif intègre.
On note $\mathcal{LA}_\mathcal{K}$ l'ensemble des paramétrisations de $E$ qui transforment tout connexe par arcs en partie connexe de $E$, et $\mathcal{LT}_\mathcal{K}$ l'ensemble des paramétrisations de $E$ qui transforment tout connexe pour la topologie usuelle en partie connexe de $E$, autrement dit : 
\[\mathcal{LA}_\mathcal{K}=\{p\in Param(E), p\in\mathbf{Cnct}((U_p,\mathcal{KA}_p) ,(E,\mathcal{K}))\}
\]
et
\[\mathcal{LT}_\mathcal{K}=\{p\in Param(E), p\in\mathbf{Cnct}((U_p,\mathcal{KT}_p) ,(E,\mathcal{K}))\}.
\]
Si $\mathcal{K}$ admet un système $\mathcal{G}$ de générateurs qui soient tous de cardinal inférieur ou égal à la puissance du continu, autrement dit s'il existe un ensemble $\mathcal{G}\subset\mathcal{K}$ tel que
\[\left\lbrace \begin{tabular}{l}
\, $\mathcal{K}=[\mathcal{G}]$, \\ 
et\\
\, $\forall G\in \mathcal{G}, card(G)\leq \mathfrak{c}$,
\end{tabular} \right.\]
alors \[U_{DC}(<\mathcal{W}_\mathcal{G}>)=U_{DC}(<\mathcal{LT}_\mathcal{K}>)=U_{DC}(<\mathcal{LA}_\mathcal{K}>)=\mathcal{K},\]
où, $\mathcal{G}$ désignant un système de générateurs ayant la propriété indiquée, on a posé 
\[\mathcal{W}_\mathcal{G}=\bigcup_{G\in \mathcal{G}}{\mathcal{W}_G}\]
avec, pour tout $G\in\mathcal{G}$,
\[\mathcal{W}_G=\{p:\mathbf{R}\rightarrow E, p(\mathbf{R})=G \,\mathrm{et}\,\forall a\in G, \overline{p^{-1}(a)}=\mathbf{R}\}.\]
\end{prop}

\noindent \textbf{Preuve}. De l'inclusion $\mathcal{KA}_p \subset \mathcal{KT}_p$, valable pour toute paramétrisation $p$ de $E$, on déduit $\mathcal{LT}_\mathcal{K}\subset \mathcal{LA}_\mathcal{K}$. Supposant à partir de maintenant l'existence d'un système de générateurs de $\mathcal{K}$ ayant les propriétés voulues, et désignant par $\mathcal{G}$ un tel système, on a en outre $\mathcal{W}_\mathcal{G}\subset\mathcal{M}_\mathcal{K}$ : en effet, pour toute paramétrisation $p\in \mathcal{W}_\mathcal{G}$, il existe $G\in\mathcal{G}\subset\mathcal{K}$ tel que pour tout connexe\footnote{Ou connexe par arcs, puisque dans $\mathbf{R}$ les deux notions sont équivalentes.} $I$ non réduit à un point et inclus dans $U_p=\mathbf{R}$, on ait $p(I)=G$. Par conséquent, le foncteur $U_{DC}$ étant trivialement croissant\footnote{Plus il y a de paramétrisations, plus il y a de connexes par arcs.}, on a
\[
U_{DC}(<\mathcal{W}_\mathcal{G}>)\subset U_{DC}(<\mathcal{LT}_\mathcal{K}>)\subset U_{DC}(<\mathcal{LA}_\mathcal{K}>).
\]

Il nous suffit donc,  pour établir la proposition \ref{prop trois diffeologisations} de montrer que, sous les hypothèses faites, on a  nécessairement
$\mathcal{K}\subset U_{DC}(<\mathcal{W}_\mathcal{G}>)$ et $U_{DC}({<\mathcal{LA}_\mathcal{K}>})\subset \mathcal{K}$.
Commençons par établir l'inclusion
\[\mathcal{K}\subset U_{DC}(<\mathcal{W}_\mathcal{G}>).\] Il suffit pour cela de vérifier 
que l'on a $\mathcal{G}\subset U_{DC}(<\mathcal{W}_\mathcal{G}>)$, puisqu'on aura alors  $\mathcal{K}=[\mathcal{G}]\subset [U_{DC}(<\mathcal{W}_\mathcal{G}>)]=U_{DC}(<\mathcal{W}_\mathcal{G}>) $.
 Soit donc $G\in\mathcal{G}$, avec $G$ non vide, de sorte que $0<card(G)\leq\mathfrak{c}$. L'ensemble $\mathbf{R}/\mathbf{Q}$, obtenu en quotientant les groupes additifs correspondants, ayant la puissance du continu, on en déduit l'existence d'une surjection ${\pi_G}:\mathbf{R}/\mathbf{Q}\twoheadrightarrow G$. En composant avec la surjection canonique $s:\mathbf{R}\twoheadrightarrow \mathbf{R}/\mathbf{Q}$, on obtient d'abord une application $p_G={\pi_G}\circ s : \mathbf{R} \twoheadrightarrow G$ puis, par extension du codomaine à $E$, une paramétrisation de $E$ que nous noterons encore $p_G:\mathbf{R} \rightarrow E$. Or, $p_G\in  \mathcal{W}_\mathcal{G}$, car la densité des ensembles $s^{-1}(c)$ pour tout $c\in \mathbf{R}/\mathbf{Q}$  implique  $\overline{p_G^{-1}(g)}=\mathbf{R}$ pour tout $g\in G$. \textit{A fortiori} $p_G\in  <\mathcal{W}_\mathcal{G}>$ : $p_G$ est une plaque de l'espace difféologique $(E,<\mathcal{W}_\mathcal{G}>)$, et par conséquent (\cite{PIZ:2013}, art. 5.9) $p_G$ transforme tout connexe par arcs
\footnote{
Les parties connectées des ouverts $U\subset \mathbf{R}^n$ munis de leur structure difféologique canonique sont les connexes par arcs.
} en partie connectée de $(E,<\mathcal{W}_\mathcal{G}>)$. En particulier, $p_G(\mathbf{R})=G$ est donc une partie connectée de $(E,<\mathcal{W}_\mathcal{G}>)$, autrement dit : $G\in U_{DC}(<\mathcal{W}_\mathcal{G}>)$. Et puisque ce résultat reste trivialement vérifié dans le cas où $G$ est vide, on a bien $\mathcal{K}\subset U_{DC}(<\mathcal{W}_\mathcal{G}>)$.

Pour établir la deuxième inclusion dont nous avons besoin, 
\[U_{DC}(<\mathcal{LA}_\mathcal{K}>)\subset\mathcal{K},\]
remarquons d'abord que, d'après le lemme \ref{lemme engendrement des connectes de diffeo par les chemins} appliqué à 
$\mathcal{D}=<\mathcal{LA}_\mathcal{K}>$ on a
\begin{equation}\label{inclusion de UDC dans les sigma}
U_{DC}(<\mathcal{LA}_\mathcal{K}>)= [\{\sigma([0,1]), \sigma\in stPaths(E,<\mathcal{LA}_\mathcal{K}>)\}].
\end{equation}

Or, pour tout $\sigma\in stPaths(E,<\mathcal{LA}_\mathcal{K}>)$, on peut écrire d'après le lemme \ref{lemme caracterisation chemins en diffeologie engendree}  l'ensemble $\sigma([0,1])\subset E$ sous la forme
\[
\sigma([0,1])=\bigcup_{1\leq k \leq n}{p_k(q_k(]a_k,b_k[))},
\] 
avec, pour tout $k$, $q_k$ de classe $C^\infty$ et $p_k\in \mathcal{LA}_\mathcal{K}$. Par continuité de $q_k$, et par la définition de $\mathcal{LA}_\mathcal{K}$ qui implique que $p_k$ transforme les connexes par arcs en connexes de $E$, on a ${\sigma(]a_k,b_k[)}={p_k(q_k(]a_k,b_k[))}\in \mathcal{K}$. 
En outre, d'après les propriétés des intervalles $]a_k,b_k[$, on a pour tout $k\in\{1,..., n-1\}$ :
${]a_k,b_k[\cap]a_{k+1},b_{k+1}[}\neq \emptyset$, d'où $\sigma(]a_k,b_k[)\cap\sigma(]a_{k+1},b_{k+1}[)\neq \emptyset$. Ainsi, $\sigma([0,1])$ peut-il s'écrire comme l'union d'une suite finie de connexes $\in\mathcal{K}$ telle que deux connexes successifs de cette suite soit non vide. On en déduit que $\sigma([0,1])\in\mathcal{K}$, d'où \[\{\sigma([0,1]), \sigma\in stPaths(E,<\mathcal{LA}_\mathcal{K}>)\}\subset \mathcal{K},\] de sorte que 
\begin{equation}\label{inclusion des sigma dans K}
[\{\sigma([0,1]), \sigma\in stPaths(E,<\mathcal{LA}_\mathcal{K}>)\}]\subset \mathcal{K}.
\end{equation}
\pagebreak[1]
Des relations (\ref{inclusion de UDC dans les sigma}) et (\ref{inclusion des sigma dans K}), on déduit $U_{DC}(<\mathcal{LA}_\mathcal{K}>)\subset\mathcal{K}$, 
ce qui achève la  démonstration. \nopagebreak[1]
\begin{flushright}$\square$\end{flushright} 



Des propositions \ref{prop conditions necessaires de diffeologisabilite} et \ref{prop trois diffeologisations}, on déduit immédiatement le théorème suivant :

\begin{thm}\label{thm de diffeologisabilite}  Un espace connectif $(E,\mathcal{K})$  est difféologisable si et seulement s'il est intègre et qu'il admet un système $\mathcal{G}\subset\mathcal{K}$ de générateurs $G\in\mathcal{G}$ qui soient tous de cardinal $card(G)\leq \mathfrak{c}$. En particulier, tout espace connectif intègre $(E,\mathcal{K})$ de support $E$ tel que $card(E)\leq\mathfrak{c}$ est difféologisable.
\end{thm}

\begin{exm} Tout espace discret intègre est difféologisable; tout espace connectif grossier est difféologisable. Les espaces connectifs usuels\footnote{C'est-à-dire ceux qui sont associés à la topologie usuelle par le foncteur d'oubli $U_T$.} $\mathbf{R}^n$ sont difféologisables, mais non pas par la difféologie usuelle sur $\mathbf{R}^n$ puisque pour celle-ci les parties connectées sont uniquement les connexes par arcs de $\mathbf{R}^n$.
Aucun espace brunnien dont le support est de cardinal strictement plus grand que $\mathfrak{c}$ n'est difféologisable, un tel espace ne vérifiant pas les conditions nécessaires de la proposition \ref{prop conditions necessaires de diffeologisabilite}.  
\end{exm}


Pour tout  espace connectif difféologisable  $(E,\mathcal{K})$  posons 
\[{\mathcal{G}_\mathcal{K}}=\{K\in\mathcal{K}, {card(K)\leq \mathfrak{c}}\}
\,\, \mathrm{et}\,\,
\mathcal{LW}_\mathcal{K}=\mathcal{W}_{\mathcal{G}_\mathcal{K}},\]
où $\mathcal{W}_{\mathcal{G}_\mathcal{K}}$ désigne la famille de paramétrisations de $E$ qui a été associée à un tel ensemble $\mathcal{G}={\mathcal{G}_\mathcal{K}}$ dans la proposition \ref{prop trois diffeologisations}. Autrement dit,
\[\mathcal{LW}_\mathcal{K}={\{f:\mathbf{R}\rightarrow E, \exists K\in\mathcal{G}_\mathcal{K}, f(\mathbf{R})=K \,\mathrm{et}\,\forall a\in K, \overline{f^{-1}(a)}=\mathbf{R}\}}.\]
Posons de plus  
\[LW({(E,\mathcal{K})})=(E,{<\mathcal{LW}_\mathcal{K}>}),\]
\[LT({(E,\mathcal{K})})=(E,{<\mathcal{LT}_\mathcal{K}>}),\]
\[LA({(E,\mathcal{K})})=(E,{<\mathcal{LA}_\mathcal{K}>}),\]
et, pour tout morphisme connectif $(f:(E,\mathcal{K})\rightarrow(E',\mathcal{K}'))\in\overrightarrow{\mathbf{Cncd}}$, 
\[LW(f)=LT(f)=LA(f)=f.\]

\begin{prop}\label{prop fonctorialite des trois diffeologisation} Les opérateurs $LW$, $LT$ et $LA$ définis ci-dessus sur les objets et sur les flèches de la catégorie $\mathbf{Cncd}$ et à valeurs dans $\mathbf{Diff}$ sont des foncteurs $\mathbf{Cncd}\rightarrow\mathbf{Diff}$.
\end{prop}
\noindent \textbf{Preuve}. Pour prouver que $LW$ est un foncteur, il suffit de vérifier que  tout morphisme connectif $(f:(E,\mathcal{K})\rightarrow(E',\mathcal{K}'))\in\overrightarrow{\mathbf{Cncd}}$ est une application lisse de $(E,\mathcal{LW}_\mathcal{K})$ dans  $(E',\mathcal{LW}_{\mathcal{K}'})$.  La vérification de ce qu'en outre $f$ est lisse de $(E,\mathcal{LT}_\mathcal{K})$ dans $(E',\mathcal{LT}_{\mathcal{K}'})$ et de $(E,\mathcal{LA}_\mathcal{K})$ dans $(E',\mathcal{LA}_{\mathcal{K}'})$ prouvera de même la fonctorialité de $LT$ et de $LA$. Or, pour toute paramétrisation $(p:\mathbf{R}\rightarrow E)\in \mathcal{LW}_\mathcal{K}$, on a $f\circ p:\mathbf{R}\rightarrow E'$ qui vérifie $f\circ p(\mathbf{R})\in \mathcal{G}_{\mathcal{K}'}$ et, pour tout $b\in f\circ p(\mathbf{R})$, ${f^{-1}(b)\neq \emptyset} \Rightarrow {\overline{(f\circ p)^{-1}(b)}=\mathbf{R}}$, de sorte que $f\circ p\in {\mathcal{LW}_{\mathcal{K}'}}$. On se trouve alors dans un cas d'application triviale du critère {1.73} de \cite{PIZ:2013} concernant les applications lisses entre espaces difféologiques dont les structures sont engendrées par des familles données de paramétrisations, à savoir le cas où pour toute paramétrisation $p$ de la première famille, $f\circ p$ appartient à la seconde famille. De même, pour toute paramétrisation $(p:\mathbf{R}^n\supset U\rightarrow E)\in {\mathcal{LT}_\mathcal{K}}$, on a $(f\circ p)\in {\mathcal{LT}_{\mathcal{K}'}}$ puisque toute partie connexe (pour la topologie usuelle) de $U$ est transformée par composition en une partie connexe de l'espace connectif ${(E',\mathcal{K}')}$. Enfin,  pour toute paramétrisation $(p:\mathbf{R}^n\supset U\rightarrow E)\in {\mathcal{LA}_\mathcal{K}}$, on a $(f\circ p)\in {\mathcal{LA}_{\mathcal{K}'}}$ puisque toute partie connexe par arcs de $U$ est transformée par composition avec $f$ en une partie connexe de l'espace connectif ${(E',\mathcal{K}')}$.
\begin{flushright}$\square$\end{flushright}



La fonctorialité des opérateurs considérés ci-dessus conduit à s'interroger sur l'existence d'adjonctions. Commençons par remarquer que $U_{DC}$ n'admet pas d'adjoint à gauche. En effet, on vérifie facilement que le foncteur $U_{DC}$ ne préserve pas les produits\footnote{Il est connu que tout foncteur qui admet un adjoint à gauche préserve toute les limites inverses, en particulier les produits.}. Considérons par exemple l'espace difféologique $\mathbf{R}$ muni de la difféologie usuelle. Alors $\mathbf{R}\times \mathbf{R}$ est le plan $\mathbf{R}^2$ muni de la difféologie usuelle\footnote{Puisqu'une plaque de cet espace produit est une paramétrisation dont les deux projections sont elles-mêmes des plaques. Sur le produit des espaces difféologiques, voir \cite{PIZ:2013}, art. 1.55.}, dont la structure connective est constituée des connexes par arcs, qui sont en particulier connexes au sens de la topologie usuelle du plan, tandis que le carré cartésien de l'espace  connectif usuel $\mathbf{R}$ admet pour connexe toute partie de $\mathbf{R}^2$ dont les deux projections sont connexes, ce qui n'implique même pas la connexité au sens de la topologie usuelle du plan\footnote{Voir \cite{Dugowson:201012}.}.

\begin{rmq} Nous avons qualifié $U_{DC}$ de \og foncteur d'oubli\fg, car il s'agit d'un foncteur fidèle entre catégories concrètes et qu'il ne retient qu'une partie des informations contenues dans une difféologie, à savoir la donnée des parties connectées. Mais il faut remarquer que, n'admettant pas d'adjoint à gauche, ce foncteur d'oubli $U_{DC}$ a un statut bien différent des foncteurs d'oubli que l'on rencontre notamment en algèbre.
\end{rmq}


\begin{prop}$LA:\mathbf{Cncd}\rightarrow \mathbf{Diff}$ est adjoint à droite de $U_{DC}:\mathbf{Diff}\rightarrow \mathbf{Cncd}$ :
\[U_{DC}\dashv LA.\]
\end{prop}

\noindent \textbf{Preuve}. Étant donnés $(X,\mathcal{D})$ un espace difféologique, $(Y,\kappa_Y)$ un espace connectif difféologisable, et $f:U_{DC}(X,\mathcal{D})\rightarrow (Y,\kappa_Y)$ une application connective, montrons que $f$ est lisse de $(X,\mathcal{D})$ dans $LA((Y,\kappa_Y))$ : pour toute plaque $(p:U_p\rightarrow X)\in\mathcal{D}$ et tout connexe par arcs $A\subset U_p$, on a $p(A)\in \mathcal{K}_\mathcal{D}=U_{DC}(\mathcal{D})$, et donc $f\circ p (A)\in \kappa_Y$, ce qui prouve que $(f\circ p)\in \mathcal{LA}_{\kappa_Y}$. Ceci établit la lisseté\footnote{\emph{Je m'aperçois à l'instant qu'à l'adjectif lisse ne correspond aucun substantif [...] Qu'il me soit permis de créer le mot \og lisseté \fg\, pour donner une idée, aux encombrés de toute nature, de ce que peut être un corps heureux}. Amélie Nothomb, \emph{Le Sabotage amoureux}, 1993.} de $f$. D'un autre coté, donnons-nous à présent une application lisse $g:(X,\mathcal{D})\rightarrow LA(Y,\kappa_Y)$. En appliquant le foncteur $U_{DC}$ à $g$, on obtient que $g=U_{DC}(g)$ est une application connective de $U_{DC}(X,\mathcal{D})$ dans $U_{DC}(LA(Y,\kappa_Y))=(Y,\kappa_Y)$. On a ainsi établi que $f\mapsto f$ est une bijection de  $\mathbf{Cncd}(U_{DC}(X,\mathcal{D}),(Y,\kappa_Y))$ sur $\mathbf{Diff}((X,\mathcal{D}),LA(Y,\kappa_Y))$, et puisque cette bijection est trivialement na-

\pagebreak[4]
\noindent turelle,  cela prouve l'adjonction annoncée.
\begin{flushright}$\square$\end{flushright} 


\begin{rmq} L'adjonction $ U_{DC}\dashv LA$  peut être comparée à celle qui a lieu entre le foncteur d'oubli de la structure topologique d'une part, et d'autre part le foncteur qui consiste à munir tout ensemble de la topologie grossière. Et, de fait, la difféologie définie par $LA$ est la plus grossière de celles qui préservent la structure connective des espaces auxquels on applique ce foncteur.
\end{rmq}


Les résultats précédents nous donnent l'occasion d'introduire la notion d'application localement connective définie sur un espace topologique et à valeurs dans un espace connectif, et de faire quelques remarques à ce sujet.

\begin{prop} Soit $(E,\mathcal{K})$ un espace connectif difféologisable, et soit $\mathcal{LA}_\mathcal{K}\subset Param(E)$ l'ensemble des paramétrisations connectives de $E$ tel que défini dans la proposition \ref{prop trois diffeologisations}. Alors 
\[\mathcal{LA}_\mathcal{K}=<\mathcal{LA}_\mathcal{K}>.\]
\end{prop}
\noindent \textbf{Preuve}.  Soit $p\in<\mathcal{LA}_\mathcal{K}>$, une plaque de l'espace difféologique     
 $(E$, ${<\mathcal{LA}_\mathcal{K}>})$. Pour tout connexe par arcs $A\subset U_p$, $p(A)$ est une partie connectée de  $(E,<\mathcal{LA}_\mathcal{K}>)$, autrement dit $p(A)\in {U_{DC}(<\mathcal{LA}_\mathcal{K}>)}$ d'où, d'après la proposition \ref{prop trois diffeologisations}, $p(A)\in{\mathcal{K}}$. On en déduit que $p$ vérifie la propriété qui caractérise les élements de $\mathcal{LA}_\mathcal{K}$, d'où l'égalité annoncée.
\begin{flushright}$\square$\end{flushright} 

\begin{rmq} 
Étant donné $(E,\mathcal{K})$ un espace connectif difféologisable, la proposition précédente entraine que l'ensemble $\mathcal{LA}_\mathcal{K}$ est une difféologie, de sorte qu'il vérifie en particulier l'axiome de localité (\cite{PIZ:2013}, art. 1.5) : pour qu'une paramétrisation $p:\mathbf{R}^n\supset U \rightarrow E$ transforme tout connexe par arcs de $U$ en connexe de $E$, il faut et il suffit qu'il existe un recouvrement ouvert de $U$ tel que la restriction de $p$ à chacun des ouverts de ce recouvrement vérifie encore la même propriété, ce qui peut d'ailleurs se vérifier directement sans difficulté. Par contre, si on remplace \og connexe par arcs\fg\, par \og connexe\fg, on obtient un énoncé qui n'est pas satisfait pour tout espace topologique. Plus précisément, étant donné $(X,\mathcal{T}_X)$ un espace topologique, disons qu'une application $f:X\rightarrow E$ est \emph{localement connective} s'il existe un recouvrement ouvert $(X_i)_{i\in I}$ de $X$  tel que, pour tout $i\in I$, $f_{\vert X_i}$ est un morphisme connectif $X_i\rightarrow E$, où $X_i$ est muni de la structure connective induite\footnote{
On vérifie facilement que la structure connective induite sur une partie $Y\subset X$ par la structure connective $U_T(\mathcal{T}_X)$ coïncide avec la structure connective $U_T(\mathcal{T}_Y)$ associée à la topologie $\mathcal{T}_Y$ induite par $\mathcal{T}_X$ sur $Y$.
}
par $U_T(\mathcal{T}_X)$. On constate alors, comme le montre le contre-exemple suivant, qu'une application localement connective n'est pas nécessairement connective : on prend pour $(X,\mathcal{T}_X)$ le sous-espace topologique du plan $\mathbf{R}^2$ induit  par la topologie usuelle sur l'ensemble $X\subset \mathbf{R}^2$ défini par
 \[X=(\bigcup_{x\in\mathbf{Q}^*}{D_x})\cup \Delta\cup\{(0,0)\},\] 
où l'on a posé $\Delta=\mathbf{R}\times\{1\}$ et, $\mathbf{Q}^*$ désignant l'ensemble des rationnels non nuls, $D_x=\{x\}\times \mathbf{R}\subset $ pour tout $x\in\mathbf{Q}^*$; pour espace connectif $(E,\mathcal{K})$, on prend $E=\{0,1\}$ muni de la structure connective discrète; et pour application $f:X\rightarrow E$, on prend celle définie par $f((0,0))=0$ et, pour tout $x\neq \{(0,0)\}$, $f(x)=1$. Cette application n'est pas connective, car $f(X)$ n'est pas connexe alors que $X$ l'est, comme on peut le vérifier ainsi : s'agissant d'un espace métrique, il suffit\footnote{\label{note renvoi exemple externe}Voir l'exemple \ref{exemple foncteur V_T} relatif au foncteur $V_T$.} de vérifier que $X$ ne peut être recouvert par deux ouverts non vides disjoints. Soient donc $A$ et $B$ deux ouverts disjoints qui recouvrent $X$, avec par exemple $A\ni (0,0)$. Alors $A\setminus \{(0,0)\}$ et $B$ sont deux ouverts disjoints qui recouvrent le connexe par arcs $(\bigcup_{x\in\mathbf{Q}^*}{D_x})\cup \Delta$, et l'on a alors $A\setminus \{(0,0)\}\neq\emptyset \Rightarrow B=\emptyset$. Donc $X$ est connexe. Maintenant, vérifions que $f$ est localement connective. Considérons pour cela $(U_i)_{i\in\{1,2\}}$ le recouvrement ouvert de $X$ défini par $U_1=\{(x,y)\in\mathbf{R}^2, y>1/4\}\cap X$ et  $U_2=\{(x,y)\in\mathbf{R}^2, y<3/4\} \cap X$. L'application $f$ étant constante sur $U_1$, elle y est nécessairement connective. La connectivité de $f$ sur $U_2$ résulte du fait que $f$ est constante sur chacune des composantes connexes de $U_2$, à savoir d'une part les segments $D_x\cap U_2$, et d'autre part le singleton  $\{0\}$. Finalement, on a bien montré que $f$ est localement connective mais n'est pas connective. Remarquons que l'espace métrique $X$ de cet exemple n'est pas localement connexe; nous laissons ici ouvert le problème de savoir si, dans le cas où l'on suppose $X$ localement connexe, en particulier si $X$ est un ouvert d'un espace $\mathbf{R}^n$, toute application localement connective $f:(X,\mathcal{T})\rightarrow (E,\mathcal{K})$ est, ou non, nécessairement connective de $U_T(X,\mathcal{T})$ dans $(E,\mathcal{K})$.
\end{rmq}

\paragraph{Remerciements.}

La rédaction de cet article a bénéficié des discussions que j'ai pu avoir avec de nombreuses personnes, et je souhaite en particulier remercier ici Anatole Khélif, Patrick Iglesias-Zemmour, Jacques Riguet, Saab Abou-Jaoudé et Andrée Ehresmann,  ainsi que les participants du séminai\-re CLE (Catégories, Logique, Etc...) dirigé à Paris VII par Anatole Khélif.



\bibliographystyle{plain}
%


%

\end{document}